\documentclass{elsart}
\usepackage{amssymb,amsbsy,amsmath,amsfonts,amssymb,amscd}
\usepackage{latexsym,euscript,epsfig}
\usepackage{pstricks}
\usepackage{array,delarray}

\newif\ifdessin
\dessintrue

\ifdessin

\usepackage{pspicture}
\usepackage[matrix,ps,xdvi]{xypic} 

\newdimen\vcadre\vcadre=0.1cm 
\newdimen\hcadre\hcadre=0.1cm 
\def\GrTeXBox#1{\vbox{\vskip\vcadre\hbox{\hskip\hcadre%
      $#1$%
   \hskip\hcadre}\vskip\vcadre}}
\def\arx#1[#2]{\ifcase#1 \relax \or%
  \ar @{-}[#2]  \or%
  \ar @2{-}[#2] \or%
  \ar @{--}[#2] \or%
  \ar @2{.}[#2] \or%
  \ar @{~}[#2]  \fi}

\fi

\def\<{\langle\,}
\def\>{\,\rangle}

\def\shuff#1#2{\mathbin{
      \hbox{\vbox{
        \hbox{\vrule
              \hskip#2
              \vrule height#1 width 0pt
               }%
        \hrule}%
             \vbox{
        \hbox{\vrule
              \hskip#2
              \vrule height#1 width 0pt
               \vrule }%
        \hrule}%
}}}
\def\shuffl{{\mathchoice{\shuff{7pt}{3.5pt}}%
                        {\shuff{6pt}{3pt}}%
                        {\shuff{4pt}{2pt}}%
                        {\shuff{3pt}{1.5pt}}}}%
\def\shuffle{\, \shuffl \,}
\def\sconc{\bullet}
 
\def\Sylv{{\rm Sylv}}   
\def\Conv{{\rm Conv}}   
\def\DeSh{{\rm DeSh}}   
\def\SG{{\mathfrak S}}  
\def\Sh{{\rm Sh}}       
\def\DC{{\rm Dec}}      
\def\maj{{\rm maj\,}}   
\def\std{{\rm Std\, }}  
\def\shape{{\rm shape\,}} 
\def\Sym{{\bf Sym}}     
\def\NCSF{{\bf Sym}}    
\def\FQSym{{\bf FQSym}} 
\def\FSym{{\bf FSym}}   
\def\QSym{{\rm QSym}}   
\def\PBT{{\bf PBT}}     
\def\sym{{\it sym}}     
\def\Prima{{\bf M}}     
\def\Primb{{\bf N}}    
\def\T{{\mathcal T}}    
\def\pp{{\mathcal P}}   
\def\qq{{\mathcal Q}}   
\def\F{{\bf F}}         
\def\G{{\bf G}}         
\def\Pp{{\bf P}}        
\def\Qq{{\bf Q}}        
\def\SS{{\bf S}}        
\def\HH{{\bf H}}        
\def\EE{{\bf E}}        
\def\qandq{\qquad\text{and}\qquad}
\def\NN{{\mathbb N}}    
\def\ZZ{{\mathbb Z}}    
\def\CC{{\mathbb C}}    
\def\KK{{\mathbb K}\, } 
\def\tensor{\otimes}
\def\Id{\operatorname{Id}}

\begin{document}

\begin{frontmatter}
\title{The Algebra of Binary Search Trees}
    
\author[moscou]{F. Hivert\thanksref{kratt}},
\author[mlv]{J.-C. Novelli\thanksref{kratt}}, and
\author[mlv]{J.-Y. Thibon\thanksref{kratt}}
\address[moscou]{
Laboratoire Franco-Russe de Math\'ematiques, Independent University of Moscow,
11, Bolchoi Vlassesky per., 121002, Moscow \\
E-mail: \texttt{Florent.Hivert@univ-mlv.fr}}
\address[mlv]{
Institut Gaspard Monge, Universit\'e de Marne-la-Vall\'ee, \\
77454 Marne-la-Vall\'ee cedex \\
E-mails: \texttt{novelli@univ-mlv.fr, jyt@univ-mlv.fr}\\
Corresponding author: Jean-Christophe Novelli}
\thanks[kratt]{This project has been partially supported by EC's IHRP
Programme, grant HPRN-CT-2001-00272, ``Algebraic Combinatorics in Europe".}

\begin{abstract}
We introduce a monoid structure on the set of binary search trees, by a
process very similar to the construction of the plactic monoid, the
Robinson-Schensted insertion being replaced by the binary search tree
insertion.
This leads to a new construction of the algebra of Planar Binary Trees of
Loday-Ronco, defining it in the same way as Non-Commutative Symmetric
Functions and Free Symmetric Functions.
We briefly explain how the main known properties of the Loday-Ronco algebra
can be described and proved with this combinatorial point of view, and then
discuss it from a representation theoretical point of view, which in turns
leads to new combinatorial properties of binary trees.
\end{abstract}
\end{frontmatter}

\tableofcontents

\section{Introduction}

There are certain analogies between the combinatorics of binary trees
and that of Young tableaux. For example, the decreasing labelings of a given
binary tree and the standard Young tableaux of a given shape are both
counted by a ``hook-length formula'' (see~\cite{FRT,St,Kn,NPS}), both
admitting natural $q$-analogs (see~\cite{St,BW1}). It is also known that one
can construct from the Robinson-Schensted correspondence a Hopf algebra
$\FSym$ whose basis is given by standard Young tableaux~\cite{PR}. Moreover, a
natural realization of this algebra by means of non-commutative polynomials
gives a enlightening proof of the Littlewood-Richardson rule
(see~\cite{Loth,DHT,NCSF6}).
In this realization, each tableau $t$ of shape $\lambda$ is interpreted as a
homogeneous polynomial of degree $n=|\lambda|$, whose commutative image is the
Schur function $s_\lambda$.

\smallskip
Recently, Loday and Ronco have introduced a Hopf algebra whose basis is the
set of planar binary trees (see~\cite{LR1,LR2}).
This algebra, denoted by $\PBT$, can as the previous one be realized
in the free associative algebra (see~\cite{DHT,NCSF6,HNT,HNT2}).
Each complete binary tree with $n$ internal nodes (or, equivalently, each
binary tree with $n$ nodes) is represented by a homogeneous polynomial of
degree $n$ in some non-commutative indeterminates.
Actually, both algebras, tableaux and binary trees, were originally defined as
Hopf sub-bialgebras of the bialgebra of permutations of \cite{MvR}, which has
been realized in~\cite{DHT,NCSF6} as the algebra of free quasi-symmetric
functions $\FQSym$.

\smallskip
This realization suggests that the algebra of Young tableaux $\FSym$ and the
algebra of planar binary trees $\PBT$ might be two particular cases of the
same construction, both based on the existence of a Robinson-Schensted-like
correspondence and a plactic-like monoid. We exhibit such a construction, the
sylvester monoid. All this process seems even more important since there is a
third example which fits in this setting: the pair of mutually dual Hopf
algebras $(\Sym,\QSym)$ (Noncommutative Symmetric Functions and
Quasi-Symmetric Functions), which corresponds to the hypoplactic monoid
(see~\cite{NCSF4,Nov}) and the Krob-Thibon correspondence.

\smallskip
Since our realization does yield a Hopf subalgebra of $\FQSym$, the basic
properties of $\PBT$ can be derived in a very natural and straigthforward
way.

The algebra of permutations $\FQSym$ is naturally equipped with a scalar
product (see~\cite{DHT,NCSF6}) 
and it is known that the integers
\begin{equation}
c_{I,J} = \< R_I, R_J \> \quad |I|=|J|=n\,,
\end{equation}
where $R_I$ stands for the noncommutative ribbon Schur function of shape $I$,
can be interpreted as Cartan invariants of the $0$-Hecke algebra: the
coefficient $c_{I,J}$ is equal to the multiplicity of the simple module $S_I$
in the indecomposable projective module $P_J$ (see~\cite{NCSF4}).

The analogy between ribbon-Schur functions and the natural basis $\Pp_T$
of $\PBT$, as presented in the sequel, allows one to wonder whether it is
possible to interpret in the same way the integers 
\begin{equation}
c_{T,U} = \< \Pp_T, \Pp_U \>
\end{equation}
This question leads us to conjecture the existence of a tower of algebras over
$\CC$: $A_0\subset A_1 \subset\cdots\subset A_n\subset\cdots$ with
$\dim_{\CC}(A_n)=n!$, whose simple modules $S_T$ (and, accordingly, the
indecomposable projective modules $P_T$) would be indexed by the  binary
trees of size $n$.
To get the complete analogy with the classical case, one should also describe
the embeddings $A_p\otimes A_q \hookrightarrow A_{p+q}$ and identify
$\PBT$ and $\PBT^*$ as the direct sums of the Grothendieck groups
\begin{equation}
{\mathcal G}=\bigoplus_{n\ge 0}G_0(A_n)\,,\quad
{\mathcal K}=\bigoplus_{n\ge 0}K_0(A_n)
\end{equation}
in such a way that the products of $\PBT$ and $\PBT^*$ written in the
appropriate basis would correspond to the induction process from
$A_p\otimes A_q$ to $A_{p+q}$.

Up to a simple reodering of the matrix, we propose a good candidate for the
matrices of Cartan invariants of a tower of algebras associated with $\PBT$
and present a precise conjecture on the induction process, suggested by some
new combinatorial properties of binary trees arising in the analysis of these
matrices.

\smallskip
This paper is organized as follows: the second Section recalls the basic
definitions and some classical combinatorial algorithms.
In Section~3, we present the algebra of planar binary trees and explain how it
leads to the definition of the sylvester monoid. We then present how one can
recover the basic properties of $\PBT$ in our setting.
In Section~4, we build new bases and derive further properties of $\PBT$.
Section~5 details the outcome of our representation theoretical 
investigation of $\PBT$. Finally, Section~7 gives some transition matrices
between different bases of $\PBT$.

The main results of this paper have been announced in~\cite{HNT,HNT2}.

\subsection*{Acknowledgments}
The computations of this paper have been done with the package MuPAD-Combinat,
part of the MuPAD project and freely available at the URL\\
{\tt http://mupad-combinat.sourceforge.net}.
The authors would like to thank Jean-Louis Loday and Mar\'\i a Ronco for
helpful discussions and comments on the first versions.
Thanks to Nicolas Thi\'ery for his help during the initial programming phase.

\section{Basic definitions and notations}

\subsection{Alphabets, words, and products}

In all the paper, we will assume that we are given a totally ordered infinite
alphabet $A$, represented either by $\{a,b,c,\ldots\}$ or by
$\{1,2,3,\ldots\}$.
Nevertheless, in some examples or for some constructions, a finite alphabet is
needed. This does not change the formulas, the only difference being that some
terms vanish in an obvious way.

\smallskip
The free associative algebra over an alphabet $A$, \emph{i.e.}, the algebra
spanned by words, the product being the concatenation, is denoted by
$\KK\<A\>$, and its unity is denoted by $\epsilon$. Here, $\KK$ is some field
of characteristic zero.

A \emph{permutation} is a word without repetition on an initial interval of
the alphabet. We shall also make use of a modification of the concatenation
product, so that, starting from two words that are permutations, one gets a
permutation.
For a word $w=x_1x_2\cdots x_n$ on the alphabet $\{1,2,\ldots\}$ and an
integer $k$, denote by $w[k]$ the word $(x_1+k)(x_2+k)\cdots (x_n+k)$, as
\emph{e.g.}, $312[4]=756$.
The \emph{shifted concatenation} of two words $u$ and $v$ is defined as
\begin{equation}
u\bullet v= u\cdot (v[k])
\end{equation}
where $k$ is the length of $u$.

\smallskip
There is another algebraic structure on $\KK\<A\>$ known as the
\emph{shuffle product}.
Let $w_1$ and $w_2$ be two words. Then the shuffle $w_1 \shuffle w_2$
is recursively defined by
\begin{itemize}
\item $w_1 \shuffle \epsilon  = w_1$,\quad $\epsilon \shuffle w_2  = w_2$,
\item $au \shuffle bv = a(u \shuffle bv) + b(au \shuffle v)$,  
\end{itemize}
where $a$, $b$ are letters, and $u$, $v$ words.

For example,
\begin{equation}
12\shuffle 43 = 1243 + 1423 + 1432 + 4123 + 4132 + 4312\,.
\end{equation}

Note that if, as in the previous example, one shuffles a permutation and
another one, shifted  by the size of the first one, one obtains a sum of
permutations.
This process is called the \emph{shifted shuffle}.

\subsection{Standardization}

The symmetric group on $n$ letters will be denoted by $\SG_n$, and its algebra
by $\KK[\SG_n]$.

Let us recall the standardization process sending a word to a permutation.

Let $A=\{a<b<\cdots\}$ be a totally ordered infinite alphabet. With each
word $w$ of $A^*$ of length $n$, we associate a permutation $\std(w)\in\SG_n$
called the \emph{standardized} of $w$ defined as the permutation obtained by
iteratively scanning $w$ from left to right, and labelling $1,2,\ldots$ the
occurrences of its smallest letter, then numbering the occurrences of the next
one, and so on. Alternatively, $\std(w)$ is the permutation having the same
inversions as $w$.

For example, $\std({abcadbcaa})=157296834$:
\begin{equation}
\begin{array}{ccccccccc}
a  & b  & c  & a  & d  & b  & c  & a  & a\\
a_1& b_5& c_7& a_2& d_9& b_6& c_8& a_3& a_4\\
1   &5   &7   &2   &9   &6   &8   &3   &4
\end{array}
\end{equation}
\medskip

\subsection{The weak order}

The \emph{weak order} (also called \emph{right permutohedron order}) is the
order on permutations obtained by defining the successors of a permutation
$\sigma$ as the permutations $\sigma\cdot s_i$ if this permutation has more
inversions than $\sigma$, where $s_i=(i,i+1)$ exchanges the numbers at places
$i$ and $i+1$ of $\sigma$ (see Figure~\ref{permu3}).

\begin{figure}[ht]
\centerline{
\newdimen\vcadre\vcadre=0.2cm 
\newdimen\hcadre\hcadre=0.2cm 
\setlength\unitlength{1.2mm}
$\xymatrix@R=.7cm@C=5mm{
      & *{\GrTeXBox{123}}\arx1[dl]\arx2[dr] &  \\
*{\GrTeXBox{213}}\arx2[d] &  & *{\GrTeXBox{132}}\arx1[d] \\
*{\GrTeXBox{231}}\arx1[dr] &  & *{\GrTeXBox{312}}\arx2[dl] \\
      & *{\GrTeXBox{321}} &  \\
}$
}
\caption{\label{permu3}The weak order of $\SG_3$.}
\end{figure}

\subsection{Permutations and saillances}

The \emph{saillances} of a permutation $\sigma$ of size $n$ are the $i\leq n$
such that all the elements to the right of $i$ in $\sigma$ are smaller than
$i$.
For example, the saillances of $893175624$ are, read from right to left, $4$,
$6$, $7$, and $9$.

For technical reasons, we build the \emph{saillance sequence} associated with
the saillances of a permutation by recording the positions of the saillances
in \emph{decreasing} order.
So the saillance sequence of $893175624$ is $(9,7,5,2)$.

\subsection{Transition matrices}

As we shall need many different bases of $\PBT$ and of various algebras, we
introduce a notation for the matrices expressing one basis into another.
The matrix $M_{A,B}$ is the matrix whose $i$th column expresses the $i$th
element of the basis $A$ as linear combination of elements of the basis $B$.

\section{The Algebra of Planar Binary Trees}

In all the paper, a \emph{planar binary tree} (or binary tree) is an
incomplete planar rooted binary tree: a binary tree is either void
($\emptyset$) or a pair of possibly void binary trees grafted on an
\emph{internal node}. The size of a tree is the number of its nodes.
The number of planar binary trees of size $n$ is the Catalan number
\begin{equation}
C_n := \frac{\binom{2n}{n}}{n+1}\,.
\end{equation}

A \emph{labeled tree} is a tree with a label attached to each node, the label
being taken either in the alphabet $A$ or in $\NN$.

A \emph{right comb tree} (resp. \emph{left comb tree}) is a tree having a
sequence of right (resp. left) edges starting from the root and trees with
only left (resp. right) edges attached to the previous nodes.

\subsection{The Loday-Ronco Algebra}

In~\cite{LR1}, Loday and Ronco introduced a Hopf Algebra of planar binary
trees, arising in their study of  dendriform algebras. Actually, this algebra
is the free dendriform algebra over one generator. In the same paper, they
proved that it is a subalgebra of the convolution algebra of permutations,
studied by Reutenauer~\cite{Reut}, Malvenuto-Reutenauer~\cite{MvR}, and
Poirier-Reutenauer \cite{PR}. We will first present this algebra in our
setting and then get back to our construction of the algebra of planar binary
trees. It will be denoted by $\PBT$, standing for \emph{Planar Binary Trees}.

\subsection{Free quasi-symmetric functions}

In~\cite{MvR}, Malvenuto and Reutenauer made a combinatorial study of the
\emph{convolution} of permutations, defined from the interpretation of
permutations as elements of the endomorphism algebra of the bialgebra
$\KK[A^*]$. 
In particular, they explicited the coproduct which endows the space
$\KK[\SG]:=\oplus_{n\geq0}\KK[\SG_n]$ (where $\SG_n$ denotes the symmetric
group) with a Hopf algebra structure.

Their theory can be significantly simplified if one embeds $\KK[\SG]$ in
$\KK\<A\>$ as in~\cite{NCSF6}.
This also sheds some light on the connection between this algebra and the
algebra of Quasi-Symmetric Functions defined by Gessel in~\cite{Ges}.

The image of $\KK[\SG]$ under this embedding is called the algebra of
\emph{Free Quasi-Symmetric functions} over $A$ and denoted by $\FQSym(A)$ or
simply by $\FQSym$ if there is no ambiguity.
Its natural basis $\F_\sigma$, where $\sigma$ runs over all permutations, is
given by the following construction.

\smallskip
\begin{defn}
Let $\sigma$ be a permutation.
The \emph{Free Quasi-Ribbon} $\F_\sigma$ is the noncommutative polynomial
\begin{equation}
\F_\sigma := \sum_{w;\,\,\std(w)=\sigma^{-1}} w\,,
\end{equation}
where $\std(w)$ denotes the standardized word of $w$ and $w$ runs over the
words on~$A$.
\end{defn}

For example, on the alphabet $\{1,2,3\}$,
\begin{equation}
\F_{123}= 111 + 112 + 113 + 122 + 123 + 133 + 222 + 223 + 233 + 333\,,
\end{equation}
\begin{equation}
\F_{213} = 212 + 213 + 313 + 323\,,
\end{equation}
\begin{equation}
\F_{312} = 221 + 231 + 331 + 332\,.
\end{equation}

With the help of the shifted shuffle, one easily describes the product of free
quasi-ribbon functions.

\begin{prop}
The free quasi-ribbons span a $\ZZ$-subalgebra of the free associative
algebra. Their product is given by the following formula.
Let $\alpha\in\SG_k$ and $\beta\in\SG_l$. Then
\begin{equation}
\label{MULTF}
\F_\alpha \F_\beta = \sum_{\sigma\in \alpha\shuffle(\beta[k])}\F_\sigma\,.
\end{equation}
\end{prop}

This algebra is in fact a Hopf algebra, the coproduct being defined as
follows. Let $A'$ and $A''$ be two mutually commuting ordered alphabets.
Identifying $F\otimes G$ with $F(A')G(A'')$, we set $\Delta(F)=F(A'\oplus
A'')$, where $\oplus$ denotes the ordered sum. Clearly, this is an algebra
homomorphism and thus defines a coproduct compatible with the product.

\begin{prop}
The coproduct of $\F_\sigma$ is given by
\begin{equation}
\label{CopF}
\Delta\F_{\sigma} = \sum_{u\cdot v=\sigma}{\F_{\std(u)} \otimes
\F_{\std(v)}}\,.
\end{equation}
\end{prop}
Moreover, $\FQSym$ is a self-dual Hopf algebra. One can see it by setting
$\G_\sigma=\F_{\sigma^{-1}}$ as a basis of $\FQSym^*$ and defining the scalar
product by
\begin{equation}
\<\F_\sigma,\G_\tau\>=\delta_{\sigma,\tau}\,.
\end{equation}

Let us recall that the convolution algebra of permutations  is the 
graded dual
$\FQSym^*$: the product of $\G$ functions is the so-called
\emph{convolution} of permutations. It consists in taking the inverses of both
permutations, make their shifted shuffle and invert the resulting  permutations.

\subsection{The sylvester monoid}

Let us now get back to the Loday-Ronco algebra.
In~\cite{LR1}, Loday and Ronco proved that $\PBT$ is a subalgebra of
the convolution algebra and gave an explicit embedding via the construction
of the decreasing tree of a permutation as described, for example,
in~\cite{Kn,Stan1}.

\begin{defn}
\label{ArbreDecr}
A \emph{decreasing tree} $T$ is a labeled tree such that the label of each
internal node is greater than the label of its children. If the labels are the
integers from 1 to  the number $n$ of nodes of $T$, 
we say that $T$ is a \emph{standard decreasing tree}.
\end{defn}

\begin{defn}
\label{w2ArbreDecr}
Let $w$ be a word without repetition. Its \emph{decreasing tree $\T(w)$} is
obtained as follows: the root is labeled by the greatest letter, $n$ of~$w$,
and if $w=unv$, where $u$ and $v$ are words without repetition, the left
subtree is $\T(u)$ and the right subtree is $\T(v)$.
\end{defn}

\begin{note}{\rm
\label{Infdec}
The left infix reading (recursively read the left subtree, the root and the
right subtree) of $\T(w)$ is $w$.
}
\end{note}

For example, the decreasing tree of $25481376$ is
\ifdessin
\begin{equation}
\label{dt254}
\T(25481376) =
  \vcenter{\tiny\xymatrix@C=1mm@R=1mm{
  &    &&  *=<10pt>[o][F-]{8}\ar@{-}[dll]\ar@{-}[drr] \\
  &  *=<10pt>[o][F-]{5}\ar@{-}[dl]\ar@{-}[dr] && &&
  *=<10pt>[o][F-]{7}\ar@{-}[dl]\ar@{-}[dr] \\
  *=<10pt>[o][F-]{2} &&  *=<10pt>[o][F-]{4}  && *=<10pt>[o][F-]{3} \ar@{-}[dl]
& & *=<10pt>[o][F-]{6} \\
         &    &             & *=<10pt>[o][F-]{1}
         }}\,.
\end{equation}
\fi

Loday and Ronco defined an embedding of $\PBT$ in the Malvenuto-Reutenauer
Hopf algebra by expressing a basis of $\PBT$ that we shall denote by $\Pp_T$,
as
\begin{equation}
\label{LodayP2F}
\Pp_T := \sum_{\sigma ;\,\,\shape(\T(\sigma))=T} \sigma\,.
\end{equation}
where $T$ is a non-labeled tree and $\shape(T)$ is the shape of the tree $T$.

This gives an embedding in $\FQSym$, which reads
\begin{equation}
\label{NousP2F}
\Pp_T = \sum_{w;\,\,\shape(\T(w))=T} w = \sum_{\sigma ; \shape(\pp(\sigma))=T}
\F_\sigma\,,
\end{equation}
where $\pp$ is a simple algorithm: it is the well-known \emph{binary search
tree insertion}, such as presented, for example, by Knuth in~\cite{Kn}.
Lemma~\ref{ppTformes} will show that this definition is consistent.

\begin{defn}
\label{ArbreBin}
A \emph{right strict binary search tree} $T$ is a labeled binary tree such
that for each internal node $n$, its label is greater than or equal to the
labels of its left subtree and strictly smaller than the labels of its right
subtree.
\end{defn}

\begin{defn}
\label{w2ArbreBin}
Let $w$ be a word. Its \emph{binary search tree $\pp(w)$} is obtained as
follows: reading $w$ from \emph{right to left}, one inserts each letter in a
binary search tree in the following way: if the tree is empty, one creates a
node labeled by the letter ; otherwise, this letter is recursively inserted in
the left (resp. right) subtree if it is smaller than or equal to (resp.
strictly greater than) the root.
\end{defn}

\begin{note}
\label{Infix}{\rm
The left infix reading of $\pp(w)$ is the non-decreasing permutation of $w$.
This is the well-known algorithm of binary search tree sorting.
}
\end{note}

Figure~\ref{fig.insert} shows the binary search tree of $cadbaedb$.

\ifdessin
\begin{figure}[ht]
\begin{multline*}
\begin{CD}
  \emptyset
  @>b>>
  \vcenter{\tiny\xymatrix@C=1mm@R=1mm{*=<10pt>[o][F-]{b}}}
  @>d>>
  \vcenter{\tiny\xymatrix@C=1mm@R=1mm{ *=<10pt>[o][F-]{b}\ar@{-}[dr] \\ &
      *=<10pt>[o][F-]{d} }}
  @>e>> \vcenter{\tiny\xymatrix@C=1mm@R=1mm{
      *=<10pt>[o][F-]{b}   \ar@{-}[dr] \\
      & *=<10pt>[o][F-]{d} \ar@{-}[dr] \\
      & & *=<10pt>[o][F-]{e} \\
    }}
    @>a>> {}
  \vcenter{\tiny\xymatrix@C=1mm@R=1mm{
      & *=<10pt>[o][F-]{b}  \ar@{-}[dl] \ar@{-}[dr] \\
      *=<10pt>[o][F-]{a} & & *=<10pt>[o][F-]{d} \ar@{-}[dr] \\
      & & & *=<10pt>[o][F-]{e} \\
  }}
\end{CD}
\\
\begin{CD}
   {} @>b>>
  \vcenter{\tiny\xymatrix@C=1mm@R=1mm{
      & *=<10pt>[o][F-]{b}  \ar@{-}[dl] \ar@{-}[dr] \\
      *=<10pt>[o][F-]{a} \ar@{-}[dr] & & *=<10pt>[o][F-]{d} \ar@{-}[dr] \\
      & *=<10pt>[o][F-]{b} & & *=<10pt>[o][F-]{e} \\
    }} @>d>>
  \vcenter{\tiny\xymatrix@C=0mm@R=1mm{
      &    &&  *=<10pt>[o][F-]{b}\ar@{-}[dll]\ar@{-}[drr] \\
      & *=<10pt>[o][F-]{a}\ar@{-}[dr] && &&
      *=<10pt>[o][F-]{d}\ar@{-}[dl]\ar@{-}[dr] \\
      &&  *=<10pt>[o][F-]{b} && *=<10pt>[o][F-]{d} && *=<10pt>[o][F-]{e} \\
    }} @>a>>
  \vcenter{\tiny\xymatrix@C=0mm@R=1mm{
      &    &&  *=<10pt>[o][F-]{b}\ar@{-}[dll]\ar@{-}[drr] \\
      & *=<10pt>[o][F-]{a}\ar@{-}[dl]\ar@{-}[dr] && &&
      *=<10pt>[o][F-]{d}\ar@{-}[dl]\ar@{-}[dr] \\
      *=<10pt>[o][F-]{a} &&  *=<10pt>[o][F-]{b} && *=<10pt>[o][F-]{d} &&
      *=<10pt >[o][F-]{e} \\
    }}
\end{CD}
\\
\begin{CD}
  {} @>c>>
  \vcenter{\tiny\xymatrix@C=0mm@R=1mm{
      &    &&  *=<10pt>[o][F-]{b}\ar@{-}[dll]\ar@{-}[drr] \\
      & *=<10pt>[o][F-]{a}\ar@{-}[dl]\ar@{-}[dr] && &&
      *=<10pt>[o][F-]{d}\ar@{-}[dl]\ar@{-}[dr] \\
      *=<10pt>[o][F-]{a} &&  *=<10pt>[o][F-]{b} && *=<10pt>[o][F-]{d}
      \ar@{-}[dl] && *=<10pt>[o][F-]{e} \\
      & & & *=<10pt>[o][F-]{c} }}\,.
\end{CD}       
\end{multline*}
\caption{\label{fig.insert}The Binary Search Tree of $cadbaedb$.}
\end{figure}
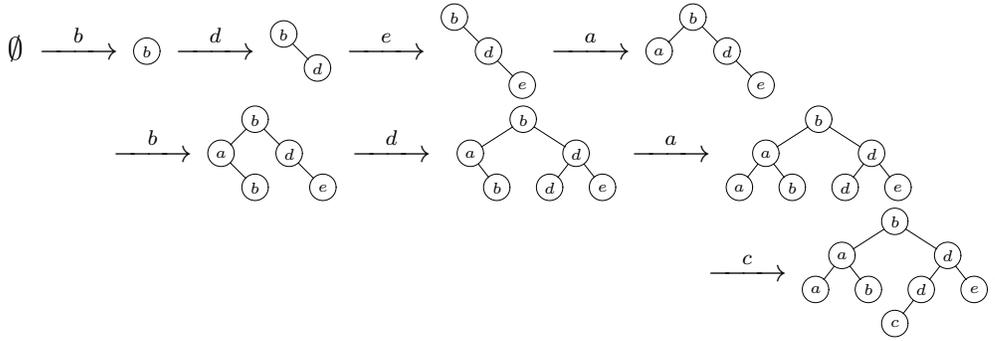
\fi

Some known examples of Hopf subalgebras of $\FQSym$ (Free Symmetric Functions
and Noncommutative Symmetric Functions) suggest to look for a monoid structure
on words to explain the previous construction.

\begin{defn}
Let $w_1$ and $w_2$ be two words. One says that they are
\emph{sylvester-adjacent} if there exists three words $u,v,w$ and three
letters $a\le b < c$ such that
\begin{equation}
w_1 = u\,ac\,v\,b\,w\qandq w_2 = u\,ca\,v\,b\,w\,.
\end{equation}

The \emph{sylvester congruence} is the transitive closure of
the relation of sylvester adjacence. 
That is, two words $u,v$ are sylvester-congruent if
there exists a chain of words
\begin{equation}
u=w_1, w_2, \dots, w_k=v\,,
\end{equation}
such that $w_i$ and $w_{i+1}$ are sylvester-adjacent for a all
$i$. In this case, we write $u \equiv_{sylv} v$.
\end{defn}
It is plain that this is actually a congruence on $A^*$.

\begin{defn}
The \emph{sylvester monoid} $\Sylv(A)$ is the quotient of the free monoid
$A^*$ by the sylvester congruence: $\Sylv(A) := A^*/\equiv_{sylv}$.
\end{defn}

For example, the classes of $w=21354$ and $w=614723$ respectively are
the sets
\begin{equation}
\label{sylvex}
\Sylv(21354) = \{21354, 21534, 25134, 52134 \}\,,
\end{equation}
\begin{equation}
\begin{split}
\Sylv(614723) = \{&126473, 162473, 164273, 164723, 612473,\\
& 614273, 614723, 641273, 641723, 647123\}\,.
\end{split}
\end{equation}

Recall that the \emph{right to left postfix reading} of a tree $T$ is the
word $w_T$ obtained by reading the right subtree, then the left and finally
the root. 
Notice that the insertion of $w_T$ with the binary search tree insertion
algorithm gives $T$ back. A word which is the postfix reading of a tree
is called the \emph{canonical word} of its sylvester class.

\begin{note}
\label{Uniq}{\rm
Thanks to Note~\ref{Infix}, one can easily see that there only is one binary
search tree of a given shape labeled by a permutation. In the sequel,
this element is called the \emph{standard canonical element of the tree}.
It now makes sense to define the canonical element of an unlabeled tree as the
canonical element of its unique binary search tree labeled by a permutation.
}
\end{note}

\begin{thm}
\label{sylv-dt}
Let $w_1$ be a word. Then $\pp(w_1)=T$ if and only if $w_1$ and $w_T$ are
sylvester-congruent.
\end{thm}

\begin{pf}
If $w_1$ and $w_T$ are sylvester-congruent, they give the same result by the
binary search tree algorithm since it is always the case for two words $w_1$
and $w_2$ which are sylvester-adjacent.

\smallskip
Conversely, any word $w$ is sylvester-congruent to a word $w_T$. Indeed, by
induction, one can assume that $w=a\cdot w_1$ where $a$ is a letter and $w_1$
a canonical word. It is then easy to see that if $a$ is smaller than or equal
to the last letter of $w_1$, it can move to the left sub-tree of $w_T$. The
result follows by induction.
\qed
\end{pf}

To get our version of the Loday-Ronco algebra, there only remains
to prove that Formula~(\ref{NousP2F}) holds.
Indeed, if we prove
\begin{equation}
\Pp_T = \sum_{\sigma;\,\, \shape(\pp(\sigma))=T} \G_{\sigma^{-1}}\,,
\end{equation}
we are done since $\FQSym$ is isomorphic to its dual, through the
identification $\F_\sigma = \G_{\sigma^{-1}}$.
This amounts to to the following Lemma.

\begin{lem}
\label{ppTformes}
Let $w$ be a word and $\sigma$ its standardized word. Then $\pp(w)$ has the
same shape as $\pp(\sigma)$ and $\T(\sigma^{-1})$.
\end{lem}

\begin{pf}
It is obvious that $\pp(w)$ and $\pp(\sigma)$ have the same shape. Let us now
prove the second part of the lemma. Let $n$ be the size of $\sigma$.
Contemplating the outputs of both insertion algorithms, the first observation
is that the size of the left subtree of $\pp(\sigma)$ is the same as the size
of the left subtree of $\T(\sigma^{-1})$: it is respectively  the number of
elements smaller than $\sigma(n)$ and the number of elements to the left of
$n$ in $\sigma^{-1}$.  Now, the proof follows by induction, since the inverse
of the standardized word of the restriction of $\sigma$ to elements smaller
than $\sigma(n)$ is the standardized word of the elements to the left of $n$
in $\sigma^{-1}$.
\qed
\end{pf}

We have now proved that the Loday-Ronco algebra of planar binary trees is
isomorphic to the algebra of sums of free quasi-ribbons over sylvester
classes. Within our description, one has
\begin{equation}
\Pp_T = \sum_{\sigma;\,\, \shape(\pp(\sigma))=T} \F_{\sigma}\,.
\end{equation}

For example, according to Equation~(\ref{sylvex}), one has:
\begin{equation}
\Pp_{\begin{picture}(5,5)
\put(1,2){\circle*{0.7}}
\put(2,1){\circle*{0.7}}
\put(1,2){\Line(1,-1)}
\put(3,3){\circle*{0.7}}
\put(3,3){\Line(-2,-1)}
\put(4,4){\circle*{0.7}}
\put(5,3){\circle*{0.7}}
\put(4,4){\Line(-1,-1)}
\put(4,4){\Line(1,-1)}
\put(4,4){\circle*{1}}
\end{picture}}
=
\Pp_{52134} = \F_{21354} + \F_{21534} + \F_{25134} + \F_{52134}\,.
\end{equation}

We have therefore realized $\PBT$ as a Hopf subalgebra of $\FQSym$ in the
same way as other interesting algebras that we will briefly recall, before
adaptating the constructions to the sylvester case.

\subsection{Analogous constructions}

The algebra of Free Symmetric Functions has been designed to give a simple and
transparent proof of the Littlewood-Richardson rule, a famous combinatorial
rule for computing tensor products of group representations which was stated
without proof in 1934, and of which no complete proof had been known until the
end of the seventies (one can find a detailed version of all this
in~\cite{NCSF6,Loth}). The monoid that will play the role of the
sylvester monoid is the well-known \emph{plactic monoid} defined by Lascoux
and Sch\"utzenberger (see~\cite{LS}) from Knuth's rewriting rules
(see~\cite{Knu}).

For later reference, let us recall that the plactic equivalence is the
congruence generated by the relations
\begin{equation}
  \left\{ \,
    \begin{array}{c@{\ }c}
      acb \, \equiv  \, cab, & \quad \hbox{for \ $a\leq b<c$}, \\
      bca \, \equiv  \, bac, & \quad \hbox{for \ $a<b\leq c$}.
    \end{array}
  \right.
\end{equation}

%
Let $t$ be a standard tableau of shape $\lambda$. Let
\begin{equation}
\SS_t := \sum_{P(\sigma)=t} \F_\sigma=\sum_{Q(w)=t }w \,,
\end{equation}
where $w\mapsto (P(w),Q(w))$ is the usual Robinson-Schensted map, sending
a word to a pair of Young Tableaux of the same shape, the second one being
standard.
As pointed out in~\cite{Loth}, Sch\"utzenberger's version of the
multiplication of Schur functions, the  \emph{Littlewood-\-Richardson
rule} is equivalent to the following statement,
which shows in particular that
the free Schur functions span a subalgebra of $\FQSym$. It is called the
algebra of \emph{Free Symmetric Functions} and denoted by $\FSym$. It
provides a realization of the algebra of tableaux introduced by Poirier and
Reutenauer~\cite{PR} as a subalgebra of the free associative algebra.

\begin{prop}
\label{LRS}
Let $t'$, $t''$ be standard tableaux, and let $k$ be the number of
cells of $t'$. Then,
\begin{equation}
\SS_{t'}\SS_{t''} = \sum_{t\in \Sh(t',t'')} \SS_t\,,
\end{equation}
where $\Sh(t',t'')$ is the set of standard tableaux in the shuffle
of $t'$ (regarded as a word via its row reading) with the plactic
class of $t''[k]$.
\end{prop}
The proof of this statement is relatively easy. It follows from
simple combinatorial properties of independent interest. In fact, the
only non trivial element of the proof is the idea of the Robinson-Schensted
correspondence. If one is willing to accept it as natural, then,
the fact that the commutative image of $\SS_t$ is $s_\lambda$
can be considered as a one-line proof of the Littlewood-Richardson rule.

Notice that one can perform the same
construction, replacing  the plactic
monoid by the \emph{hypoplactic monoid} (see~\cite{NCSF4,Nov}). One then
obtains the algebra of Noncommutative Symmetric Functions (see~\cite{NCSF7}).

\subsection{A Schensted-like algorithm}

We have already mentioned 
the Robinson-Schensted correspondence, a bijection between words and pairs of 
tableaux computed by
the Schensted algorithm (see~\cite{Sch}). The same construction
generalizes to pairs composed of ribbons and quasi-ribbons in the case of the
hypoplactic monoid (see~\cite{NCSF4,Nov}). This construction also generalizes
to the sylvester case.

Let $w$ be a word. The \emph{Sylvester Schensted Algorithm} $SSA$ sends it to
the pair composed of its binary search tree $\pp(w)$ and the decreasing tree
of the inverse of its standardized $\qq(w)=\T((\std(w))^{-1})$.

For example, according to Formula~(\ref{dt254}) and Figure~\ref{fig.insert},
one has

\ifdessin
\begin{equation}
SSA(cadbaedb) = \left(
  \vcenter{\tiny\xymatrix@C=1mm@R=1mm{
  &    &&  *=<10pt>[o][F-]{b}\ar@{-}[dll]\ar@{-}[drr] \\
  &  *=<10pt>[o][F-]{a}\ar@{-}[dl]\ar@{-}[dr] && &&
  *=<10pt>[o][F-]{d}\ar@{-}[dl]\ar@{-}[dr] \\
  *=<10pt>[o][F-]{a} &&  *=<10pt>[o][F-]{b} && *=<10pt>[o][F-]{d} \ar@{-}[dl]
&& *=<10pt>[o][F-]{e} \\
         &    &             & *=<10pt>[o][F-]{c}
         }}
\ ,\quad
  \vcenter{\tiny\xymatrix@C=1mm@R=1mm{
  &    &&  *=<10pt>[o][F-]{8}\ar@{-}[dll]\ar@{-}[drr] \\
  &  *=<10pt>[o][F-]{5}\ar@{-}[dl]\ar@{-}[dr] && &&
  *=<10pt>[o][F-]{7}\ar@{-}[dl]\ar@{-}[dr] \\
  *=<10pt>[o][F-]{2} &&  *=<10pt>[o][F-]{4}  && *=<10pt>[o][F-]{3} \ar@{-}[dl]
& & *=<10pt>[o][F-]{6} \\
         &    &             & *=<10pt>[o][F-]{1}
         }}
  \right)\,.
\end{equation}
\fi

Consider a pair composed of a binary search tree $t_1$ and a standard
decreasing tree $t_2$ of the same shape.
Algorithm $SSB$ sends this pair to the word obtained by reading the labels of
$t_1$ in the order of the corresponding labels in $t_2$.

\begin{thm}
\label{Schens}
The algorithm $SSA$ yields a bijection between the set of words of size $n$
and pairs composed of a binary search tree of size $n$ and a standard
decreasing tree of the same shape. The reciprocal bijection is computed by
Algorithm $SSB$.
\end{thm}

\begin{pf}
First, the  output of Algorithm~$SSA$ is a pair of trees of the same
shape, thanks to Lemma~\ref{ppTformes}.
So it maps $w$ to a pair of the right form. Moreover, thanks to
Notes~\ref{Infdec} and~\ref{Infix} and since $w$ is equal to the word obtained
by permuting its non-decreasing word with the infix reading of
$\T(\std(w)^{-1})$, one can conclude.
\qed
\end{pf}

\begin{note}
\label{sylv2poset}{\rm
If one regards a tree as a partially ordered set, the leafs being the smallest
elements and the root the greatest, the sylvester class associated with a
given tree consists in the linear extensions of this partial order.
}
\end{note}

\subsection{Basic properties of the sylvester monoid and of its algebra}

\subsubsection{Sylvester classes and the permutohedron}

First, let us describe the structure of the sylvester equivalence classes seen
as parts of the permutohedron. The first two properties are obvious thanks to
Note~\ref{sylv2poset}.

\begin{prop}
The greatest word for the lexicographic order of the sylvester class of $T$ is
$w_T$.
\end{prop}

\begin{prop}
The smallest word for the lexicographic order of the sylvester class 
indexed by a tree $T$ is the left to right postfix reading of $T$.
\end{prop}

The following property will later be useful to simplify the combinatorial
description of $\PBT$. We first recall the definition of pattern avoidance.

\begin{def}
\label{def:patav}
A permutation $\sigma$ of $\SG_n$ \emph{avoids the pattern} $\pi$ of $\SG_k$
if and only if there is no subsequence $i_{\pi(1)} < i_{\pi(2)} < \ldots <
i_{\pi(k)}$ in $[1,n]$ such that $\sigma(i_1)<\sigma(i_2)<\ldots<\sigma(i_k)$.
\end{def}

\begin{prop}
\label{motif132}
A permutation is a canonical sylvester permutation if and only if it avoids
the pattern $132$.
\end{prop}

\begin{pf}
If $\sigma$ is a canonical sylvester permutation, since it is the right to
left postfix reading of a binary search tree, it necessarily avoids the
pattern $132$.
Since the number of permutations avoiding $132$ of size $n$ is the Catalan
number $C_n$, there are as many permutations of size $n$ avoiding $132$ as
binary trees with $n$ nodes.
\qed
\end{pf}

If $\sigma$ avoids a given pattern $\pi$ then $\sigma^{-1}$ avoids the pattern
$\pi^{-1}$. Consequently:

\begin{cor}
\label{invcanon}
Let $\sigma$ be a canonical sylvester permutation. Then $\sigma^{-1}$ also is
a canonical sylvester permutation.
\end{cor}

The next property can be used to prove that the $\Pp_T$'s span a Hopf
subalgebra of $\FQSym$.
Let us examine the compatibility of the sylvester congruence with restriction
to intervals.

\begin{prop}
Assume that $u$ and $v$ are sylvester-congruent. Let $I$ be an interval of the
alphabet $A$, $I=[a_k,\dots,a_l]$.  Let $u/I$ (resp. $v/I$) be the word
obtained by erasing the letters of $u$ (resp. $v$) that are not in $I$. Then
$u/I$ and $v/I$ are sylvester congruent.
\end{prop}

\begin{pf}
This property is easily checked on the sylvester relations, which implies the
result.
\qed
\end{pf}

Now, by the very same reasoning as in~\cite{NCSF6}, prop. 3.12, using the
fact that the sylvester congruence is compatible with the restriction to
intervals and to de-standardization, one deduces that $\PBT$ is a Hopf
subalgebra of $\FQSym$. All formulas will be given in the next Section.

\begin{note}{\rm
In~\cite{LR2}, the product $\Pp_{T'}\Pp_{T''}$ is described by means of an
order on the planar binary trees, also known as the \emph{Tamari} order (see
Section~\ref{Tamari} for more details).
In our setting, this order is obtained from the weak order of the symmetric
group: it is its restriction to canonical words (see
Theorem~\ref{tamari2canon}).
It should be noticed that the results of Bj\"orner and Wachs (see~\cite{BW2})
show that the sylvester classes are intervals for the weak order.
}
\end{note}

\begin{note}{\rm
In the context of tableaux, as a simple consequence of Theorem~1
of~\cite{NCSF6}, one can define an order on tableaux as the quotient of the
weak order by Knuth's relations (this order is also considered by Bj\"orner
and Wachs).
The sets $\Sh(t',t'')$ of Proposition~\ref{LRS}, where $t'$ and $t''$ are
standard tableaux, are intervals for this order.
}
\end{note}

\subsubsection{Hook-length formula for trees}

For sake of completeness, we include Knuth's hook-length formula for trees, 
and its $q$-analog, due to Bj\"orner and Wachs.

\smallskip
The cardinality of a standard plactic class is equal to the number of
standard tableaux of a certain shape which is given by the celebrated
hook-length formula (see~\cite{FRT,NPS}). This formula admits a $q$-analog
which enumerates permutations by their major index (see~\cite{St}). In the
same way, the enumeration of the sylvester class associated with a tree $T$ of
size $n$ is given by the specialization $(q)_n\Pp_T(1,q,q^2,\cdots)$ which is
equal (see~\cite{BW1}) to
\begin{equation}
  \sum_{\pp(\sigma)=T}q^{\maj(\sigma)} =
  (q)_n\Pp_T(1,q,q^2,\cdots) =
  \frac{[n]_q!}{\prod_{\circ\in T} q^{-\delta_\circ}[h_\circ]_q}\,,
\end{equation}
where for a node $\circ$ of $T$, the coefficient $h_\circ$ is the size of the
subtree rooted at $\circ$ and $\delta_\circ$ the size of its right subtree.
Here we use the standard notations of the $q$-calculus, that is, for any
integer $n$, the $q$-integer $[n]_q$ is equal to $1+q+\dots+q^{n-1}$, the
$q$-factorial $[n]!_q$ is the product of the corresponding $q$-integers:
$[n]!_q := [1]_q[2]_1\dots[n]_q$, and finally
\begin{equation}
(q)_n := (1-q)^n [n]!_q = (1-q) (1-q^2) \dots (1-q^n)\,.
\end{equation}

\subsection{A sylvester description of $\PBT$}

We are now in position to describe the structure of $\PBT$ in a similar
way as in \cite{NCSF6,NT} for other combinatorial Hopf algebras (in the
sense of~\cite{ABS}). In particular we will give alternative formulas to
compute the product, the coproduct and the antipode of $\PBT$  using
only combinatorial properties of the sylvester classes.

\begin{thm}
\label{ProdP}
Let $T'$ and $T''$ be two binary trees. Then
\begin{equation}
\Pp_{T'}\Pp_{T''}=\sum_{T\in \Sh(T',T'')}\Pp_T\,,
\end{equation}
where $\Sh(T',T'')$ is the set of trees $T$ such that $w_T$ occur in the
shuffle product $u\shuffle v$,  $u=w_{T'}$ and $v=w_{T''}[k]$ is the
canonical word of $T''$ shifted by the size of $T'$.
\end{thm}

\begin{pf}
The product $\Pp_{T'}\Pp_{T''}$ can be expanded into $\F_\sigma$'s using the
shifted shuffle on permutations (see Formulas~\ref{NousP2F} and~\ref{MULTF}).
It can be factored as a sum of $\Pp_T$'s thanks to the compatibility of the
sylvester relation with restriction to intervals. Finally, $\Pp_T$ arises in
the product $\Pp_{T'}\Pp_{T''}$ if and only if its canonical word $w_T$
appears in the shuffle of the sylvester classes associated with $T'$ and
$T''$. The last part of the theorem then comes from the fact that if $w$ and
$w'$ are not both canonical words, there is no canonical word in their shifted
shuffle.
\qed
\end{pf}

For example,
\begin{equation}
\Pp_{4213} \Pp_{312} = \Pp_{7421356} + \Pp_{7452136} + \Pp_{7456213} +
\Pp_{7542136} + \Pp_{7546213} + \Pp_{7564213}\,,
\end{equation}
\begin{equation}
\begin{split}
\Pp_{{\begin{picture}(4,4)
\put(1,2){\circle*{0.7}}
\put(2,1){\circle*{0.7}}
\put(1,2){\Line(1,-1)}
\put(3,3){\circle*{0.7}}
\put(4,2){\circle*{0.7}}
\put(3,3){\Line(-2,-1)}
\put(3,3){\Line(1,-1)}
\put(3,3){\circle*{1}}
\end{picture}}}
\Pp_{{\begin{picture}(3,3)
\put(1,1){\circle*{0.7}}
\put(2,2){\circle*{0.7}}
\put(3,1){\circle*{0.7}}
\put(2,2){\Line(-1,-1)}
\put(2,2){\Line(1,-1)}
\put(2,2){\circle*{1}}
\end{picture}}}
\ =\ &
\Pp_{{\begin{picture}(7,7)
\put(1,3){\circle*{0.7}}
\put(2,2){\circle*{0.7}}
\put(1,3){\Line(1,-1)}
\put(3,4){\circle*{0.7}}
\put(4,3){\circle*{0.7}}
\put(3,4){\Line(-2,-1)}
\put(3,4){\Line(1,-1)}
\put(5,5){\circle*{0.7}}
\put(5,5){\Line(-2,-1)}
\put(6,6){\circle*{0.7}}
\put(7,5){\circle*{0.7}}
\put(6,6){\Line(-1,-1)}
\put(6,6){\Line(1,-1)}
\put(6,6){\circle*{1}}
\end{picture}}} \ + \ \Pp_{{\begin{picture}(7,7)
\put(1,4){\circle*{0.7}}
\put(2,3){\circle*{0.7}}
\put(1,4){\Line(1,-1)}
\put(3,5){\circle*{0.7}}
\put(4,3){\circle*{0.7}}
\put(5,4){\circle*{0.7}}
\put(5,4){\Line(-1,-1)}
\put(3,5){\Line(-2,-1)}
\put(3,5){\Line(2,-1)}
\put(6,6){\circle*{0.7}}
\put(7,5){\circle*{0.7}}
\put(6,6){\Line(-3,-1)}
\put(6,6){\Line(1,-1)}
\put(6,6){\circle*{1}}
\end{picture}}} \ + \ \Pp_{{\begin{picture}(7,7)
\put(1,4){\circle*{0.7}}
\put(2,3){\circle*{0.7}}
\put(1,4){\Line(1,-1)}
\put(3,5){\circle*{0.7}}
\put(4,2){\circle*{0.7}}
\put(5,3){\circle*{0.7}}
\put(5,3){\Line(-1,-1)}
\put(6,4){\circle*{0.7}}
\put(7,3){\circle*{0.7}}
\put(6,4){\Line(-1,-1)}
\put(6,4){\Line(1,-1)}
\put(3,5){\Line(-2,-1)}
\put(3,5){\Line(3,-1)}
\put(3,5){\circle*{1}}
\end{picture}}} \\
&+ \ \Pp_{{\begin{picture}(7,6)
\put(1,4){\circle*{0.7}}
\put(2,3){\circle*{0.7}}
\put(1,4){\Line(1,-1)}
\put(3,5){\circle*{0.7}}
\put(4,4){\circle*{0.7}}
\put(5,3){\circle*{0.7}}
\put(4,4){\Line(1,-1)}
\put(3,5){\Line(-2,-1)}
\put(3,5){\Line(1,-1)}
\put(6,6){\circle*{0.7}}
\put(7,5){\circle*{0.7}}
\put(6,6){\Line(-3,-1)}
\put(6,6){\Line(1,-1)}
\put(6,6){\circle*{1}}
\end{picture}}} \ + \ \Pp_{{\begin{picture}(7,6)
\put(1,4){\circle*{0.7}}
\put(2,3){\circle*{0.7}}
\put(1,4){\Line(1,-1)}
\put(3,5){\circle*{0.7}}
\put(4,3){\circle*{0.7}}
\put(5,2){\circle*{0.7}}
\put(4,3){\Line(1,-1)}
\put(6,4){\circle*{0.7}}
\put(7,3){\circle*{0.7}}
\put(6,4){\Line(-2,-1)}
\put(6,4){\Line(1,-1)}
\put(3,5){\Line(-2,-1)}
\put(3,5){\Line(3,-1)}
\put(3,5){\circle*{1}}
\end{picture}}} \ + \ \Pp_{{\begin{picture}(7,6)
\put(1,4){\circle*{0.7}}
\put(2,3){\circle*{0.7}}
\put(1,4){\Line(1,-1)}
\put(3,5){\circle*{0.7}}
\put(4,4){\circle*{0.7}}
\put(5,2){\circle*{0.7}}
\put(6,3){\circle*{0.7}}
\put(7,2){\circle*{0.7}}
\put(6,3){\Line(-1,-1)}
\put(6,3){\Line(1,-1)}
\put(4,4){\Line(2,-1)}
\put(3,5){\Line(-2,-1)}
\put(3,5){\Line(1,-1)}
\put(3,5){\circle*{1}}
\end{picture}}}
\end{split}
\end{equation}

\smallskip
Let us now compute the coproduct of a $\Pp_T$ in a similar way.

\begin{thm}
\label{coprodP}
Let $T$ be a tree. The coproduct of $\Pp_T$ is given by
\begin{equation}
\Delta \Pp_T = \sum_{(T',T'')\in \DC(T)} \Pp_{T'} \otimes \Pp_{T''}\,,
\end{equation}
where $\DC(T)$ is the set of pairs of trees $(T',T'')$ such that
$w_{T'}$ (resp. $w_{T''}$) are the standardized words of elements $w_1$ (resp.
$w_2$) where $w_1\cdot w_2$ is in the sylvester class of $T$.
\end{thm}

\begin{pf}
The coproduct of a $\Pp_T$ can be expanded into $\F_\sigma$'s, then
described by using the deconcatenation of permutations (see
Formulas~(\ref{NousP2F}) and~(\ref{CopF})). It can be factored into a sum of
$\Pp_T$'s thanks to the compatibility of the sylvester relation with
de-standardization. Finally, a pair $(\Pp_{T'},\Pp_{T''})$ arises in the
coproduct of $\Pp_{T}$ if and only if both canonical words $w_{T'}$ and
$w_{T''}$ appear in the deconcatenation of an element of the sylvester class
of $T$.
\qed
\end{pf}

For example,
\begin{equation}
\begin{split}
\Delta\Pp_{4213} =& \ \ \Pp_{4213}\otimes1 +
\left(\Pp_{213}+\Pp_{231}+\Pp_{321}\right)\otimes\Pp_{1} +
\left(\Pp_{12}+\Pp_{21}\right)\otimes\Pp_{12} \\
& + \Pp_{21}\otimes\Pp_{21} + \Pp_{1}\otimes\left(\Pp_{213}+\Pp_{312}\right) +
1\otimes\Pp_{4213}\,,
\end{split}
\end{equation}
\begin{equation}
%
\begin{split}
\Delta \Pp_{{\begin{picture}(4,4)
\put(1,2){\circle*{0.7}}
\put(2,1){\circle*{0.7}}
\put(1,2){\Line(1,-1)}
\put(3,3){\circle*{0.7}}
\put(4,2){\circle*{0.7}}
\put(3,3){\Line(-2,-1)}
\put(3,3){\Line(1,-1)}
\put(3,3){\circle*{1}}
\end{picture}}} \ =\ & \Pp_{{\begin{picture}(4,4)
\put(1,2){\circle*{0.7}}
\put(2,1){\circle*{0.7}}
\put(1,2){\Line(1,-1)}
\put(3,3){\circle*{0.7}}
\put(4,2){\circle*{0.7}}
\put(3,3){\Line(-2,-1)}
\put(3,3){\Line(1,-1)}
\put(3,3){\circle*{1}}
\end{picture}}}\,\otimes\,1 \ +
\left( \ \Pp_{{\begin{picture}(3,4)
\put(1,2){\circle*{0.7}}
\put(2,1){\circle*{0.7}}
\put(1,2){\Line(1,-1)}
\put(3,3){\circle*{0.7}}
\put(3,3){\Line(-2,-1)}
\put(3,3){\circle*{1}}
\end{picture}}}
\ + \ \Pp_{{\begin{picture}(3,4)
\put(1,3){\circle*{0.7}}
\put(2,1){\circle*{0.7}}
\put(3,2){\circle*{0.7}}
\put(3,2){\Line(-1,-1)}
\put(1,3){\Line(2,-1)}
\put(1,3){\circle*{1}}
\end{picture}}}
\ + \ \Pp_{{\begin{picture}(3,4)
\put(1,3){\circle*{0.7}}
\put(2,2){\circle*{0.7}}
\put(3,1){\circle*{0.7}}
\put(2,2){\Line(1,-1)}
\put(1,3){\Line(1,-1)}
\put(1,3){\circle*{1}}
\end{picture}}}\ \right)\,\otimes\, \Pp_{{\begin{picture}(1,2)
\put(1,1){\circle*{0.7}}
\put(1,1){\circle*{1}}
\end{picture}}} \\
&
+ \left(\ \Pp_{{\begin{picture}(2,3)
\put(1,1){\circle*{0.7}}
\put(2,2){\circle*{0.7}}
\put(2,2){\Line(-1,-1)}
\put(2,2){\circle*{1}}
\end{picture}}}
\ + \ \Pp_{{\begin{picture}(2,3)
\put(1,2){\circle*{0.7}}
\put(2,1){\circle*{0.7}}
\put(1,2){\Line(1,-1)}
\put(1,2){\circle*{1}}
\end{picture}}}\ \right) \otimes\, \Pp_{{\begin{picture}(2,3)
\put(1,1){\circle*{0.7}}
\put(2,2){\circle*{0.7}}
\put(2,2){\Line(-1,-1)}
\put(2,2){\circle*{1}}
\end{picture}}}
\ + \ \Pp_{{\begin{picture}(2,3)
\put(1,2){\circle*{0.7}}
\put(2,1){\circle*{0.7}}
\put(1,2){\Line(1,-1)}
\put(1,2){\circle*{1}}
\end{picture}}}\,\otimes\, \Pp_{{\begin{picture}(2,3)
\put(1,2){\circle*{0.7}}
\put(2,1){\circle*{0.7}}
\put(1,2){\Line(1,-1)}
\put(1,2){\circle*{1}}
\end{picture}}} \\
& \ + \ \Pp_{{\begin{picture}(1,2)
\put(1,1){\circle*{0.7}}
\put(1,1){\circle*{1}}
\end{picture}}}\,\otimes\,
\left(
\ \Pp_{{\begin{picture}(3,4)
\put(1,2){\circle*{0.7}}
\put(2,1){\circle*{0.7}}
\put(1,2){\Line(1,-1)}
\put(3,3){\circle*{0.7}}
\put(3,3){\Line(-2,-1)}
\put(3,3){\circle*{1}}
\end{picture}}}
\ +
\ \Pp_{{\begin{picture}(3,3)
\put(1,1){\circle*{0.7}}
\put(2,2){\circle*{0.7}}
\put(3,1){\circle*{0.7}}
\put(2,2){\Line(-1,-1)}
\put(2,2){\Line(1,-1)}
\put(2,2){\circle*{1}}
\end{picture}}}\ \right)
\ + \ 1\,\otimes\, \Pp_{{\begin{picture}(4,4)
\put(1,2){\circle*{0.7}}
\put(2,1){\circle*{0.7}}
\put(1,2){\Line(1,-1)}
\put(3,3){\circle*{0.7}}
\put(4,2){\circle*{0.7}}
\put(3,3){\Line(-2,-1)}
\put(3,3){\Line(1,-1)}
\put(3,3){\circle*{1}}
\end{picture}}}
\end{split}
\end{equation}

Since $\PBT$ is a connected graded bialgebra of finite dimension in each
component, one can define the antipode $\nu$ of $\PBT$ without ambiguity. It
then endows $\PBT$ with the structure of a graded Hopf algebra. We will
provide its formula on the dual basis of the $\Pp_T$ functions (see
Equation~(\ref{antipQ})).

%
\section{Properties of $\PBT$}
%

\subsection{Duality}

The sylvester congruence also gives a nice characterization of the dual
algebra of $\PBT$:

\begin{thm}
The (graded) dual $\PBT^*$ of the algebra of planar binary trees $\PBT$ is
isomorphic to the image of $\FQSym^*$ under the canonical projection 
\begin{equation}
\pi :\ \CC\<A\>\longrightarrow \CC[\Sylv(A)]\simeq \CC\<A\>/\equiv_{sylv}\,.
\end{equation}
The dual basis $\Qq_T$ of $\Pp_T$ is expressed as $\Qq_T=\pi(\G_\sigma)$,
where $\sigma$ is any permutation of the sylvester class associated with $T$.
\end{thm}

Notice that all this works within the realization of $\FQSym^*$ since the
sylvester monoid is compatible with the de-standardization process.
Let us now see how to compute the product and coproduct of $\Qq_T$ functions 
by means of our formalism.

\begin{thm}
\label{prodG}
Let $T'$ and $T''$ be two trees. Then, 
\begin{equation}
\Qq_{T'} \Qq_{T''} = \sum_{T\in \Conv(T',T'')}\Qq_{T}\,,
\end{equation}
where $\Conv(T',T'')$ is the set of trees that are the binary search trees of
an element of the convolution product of $w_{T'}$ by $w_{T''}$.
\end{thm}

For example,
\begin{equation}
\begin{split}
\Qq_{21} \Qq_{312} =&\ \  \Qq_{43512} + \Qq_{45123} + \Qq_{45213} +
\Qq_{52134} + \Qq_{53124}\\
& + \Qq_{53214} + \Qq_{53412} + \Qq_{54123} + \Qq_{54213} + \Qq_{54312}\,.
\end{split}
\end{equation}
\begin{equation}
\begin{split}
\Qq_{{\begin{picture}(2,3)
\put(1,2){\circle*{0.7}}
\put(2,1){\circle*{0.7}}
\put(1,2){\Line(1,-1)}
\put(1,2){\circle*{1}}
\end{picture}}}
\Qq_{{\begin{picture}(3,3)
\put(1,1){\circle*{0.7}}
\put(2,2){\circle*{0.7}}
\put(3,1){\circle*{0.7}}
\put(2,2){\Line(-1,-1)}
\put(2,2){\Line(1,-1)}
\put(2,2){\circle*{1}}
\end{picture}}}
\ = \ &
\Qq_{{\begin{picture}(5,5) 
\put(1,3){\circle*{0.7}}
\put(2,4){\circle*{0.7}}
\put(3,2){\circle*{0.7}}
\put(4,1){\circle*{0.7}}
\put(3,2){\Line(1,-1)}
\put(5,3){\circle*{0.7}}
\put(5,3){\Line(-2,-1)}
\put(2,4){\Line(-1,-1)}
\put(2,4){\Line(3,-1)}
\put(2,4){\circle*{1}}
\end{picture}}} \ + \ \Qq_{{\begin{picture}(5,4) 
\put(1,1){\circle*{0.7}}
\put(2,2){\circle*{0.7}}
\put(2,2){\Line(-1,-1)}
\put(3,3){\circle*{0.7}}
\put(4,1){\circle*{0.7}}
\put(5,2){\circle*{0.7}}
\put(5,2){\Line(-1,-1)}
\put(3,3){\Line(-1,-1)}
\put(3,3){\Line(2,-1)}
\put(3,3){\circle*{1}}
\end{picture}}} \ + \ \Qq_{{\begin{picture}(5,4) 
\put(1,2){\circle*{0.7}}
\put(2,1){\circle*{0.7}}
\put(1,2){\Line(1,-1)}
\put(3,3){\circle*{0.7}}
\put(4,1){\circle*{0.7}}
\put(5,2){\circle*{0.7}}
\put(5,2){\Line(-1,-1)}
\put(3,3){\Line(-2,-1)}
\put(3,3){\Line(2,-1)}
\put(3,3){\circle*{1}}
\end{picture}}} \ + \ \Qq_{{\begin{picture}(5,5) 
\put(1,2){\circle*{0.7}}
\put(2,1){\circle*{0.7}}
\put(1,2){\Line(1,-1)}
\put(3,3){\circle*{0.7}}
\put(3,3){\Line(-2,-1)}
\put(4,4){\circle*{0.7}}
\put(5,3){\circle*{0.7}}
\put(4,4){\Line(-1,-1)}
\put(4,4){\Line(1,-1)}
\put(4,4){\circle*{1}}
\end{picture}}}
\ + \ \Qq_{{\begin{picture}(5,5) 
\put(1,2){\circle*{0.7}}
\put(2,3){\circle*{0.7}}
\put(3,2){\circle*{0.7}}
\put(2,3){\Line(-1,-1)}
\put(2,3){\Line(1,-1)}
\put(4,4){\circle*{0.7}}
\put(5,3){\circle*{0.7}}
\put(4,4){\Line(-2,-1)}
\put(4,4){\Line(1,-1)}
\put(4,4){\circle*{1}}
\end{picture}}} \\
&+ \ \Qq_{{\begin{picture}(5,5) 
\put(1,3){\circle*{0.7}}
\put(2,2){\circle*{0.7}}
\put(3,1){\circle*{0.7}}
\put(2,2){\Line(1,-1)}
\put(1,3){\Line(1,-1)}
\put(4,4){\circle*{0.7}}
\put(5,3){\circle*{0.7}}
\put(4,4){\Line(-3,-1)}
\put(4,4){\Line(1,-1)}
\put(4,4){\circle*{1}}
\end{picture}}} \ + \ \Qq_{{\begin{picture}(5,5) 
\put(1,3){\circle*{0.7}}
\put(2,4){\circle*{0.7}}
\put(3,2){\circle*{0.7}}
\put(4,3){\circle*{0.7}}
\put(5,2){\circle*{0.7}}
\put(4,3){\Line(-1,-1)}
\put(4,3){\Line(1,-1)}
\put(2,4){\Line(-1,-1)}
\put(2,4){\Line(2,-1)}
\put(2,4){\circle*{1}}
\end{picture}}} \ + \ \Qq_{{\begin{picture}(5,4) 
\put(1,1){\circle*{0.7}}
\put(2,2){\circle*{0.7}}
\put(2,2){\Line(-1,-1)}
\put(3,3){\circle*{0.7}}
\put(4,2){\circle*{0.7}}
\put(5,1){\circle*{0.7}}
\put(4,2){\Line(1,-1)}
\put(3,3){\Line(-1,-1)}
\put(3,3){\Line(1,-1)}
\put(3,3){\circle*{1}}
\end{picture}}} \ + \ \Qq_{{\begin{picture}(5,4) 
\put(1,2){\circle*{0.7}}
\put(2,1){\circle*{0.7}}
\put(1,2){\Line(1,-1)}
\put(3,3){\circle*{0.7}}
\put(4,2){\circle*{0.7}}
\put(5,1){\circle*{0.7}}
\put(4,2){\Line(1,-1)}
\put(3,3){\Line(-2,-1)}
\put(3,3){\Line(1,-1)}
\put(3,3){\circle*{1}}
\end{picture}}} \ + \ \Qq_{{\begin{picture}(5,5) 
\put(1,3){\circle*{0.7}}
\put(2,4){\circle*{0.7}}
\put(3,3){\circle*{0.7}}
\put(4,2){\circle*{0.7}}
\put(5,1){\circle*{0.7}}
\put(4,2){\Line(1,-1)}
\put(3,3){\Line(1,-1)}
\put(2,4){\Line(-1,-1)}
\put(2,4){\Line(1,-1)}
\put(2,4){\circle*{1}}
\end{picture}}}
\end{split}
\end{equation}

\begin{thm}
\label{coprodG}
Let $T$ be a tree. Then, 
\begin{equation}
\Delta\Qq_{T} = \sum_{(T',T'')\in\DeSh(T)} \Qq_{T'}\otimes\Qq_{T''}\,,
\end{equation}
where $\DeSh(T)$ is the set of pair of trees $(T',T'')$ such that their
canonical elements are the standardized of the restrictions of the canonical
word of $T$ to all pairs of intervals $[1,i]$ and $[i+1,n]$ for
$i\in\{0,\ldots,n\}$.
\end{thm}

For example,
\begin{equation}
\begin{split}
\Delta\Qq_{645213} =&\ \ \Qq_{645213}\otimes1 + \Qq_{45213}\otimes\Qq_{1} +
\Qq_{4213}\otimes\Qq_{21} + \Qq_{213}\otimes\Qq_{312}\\
& + \Qq_{21}\otimes\Qq_{4231} +
\Qq_{1}\otimes\Qq_{53412} + 1\otimes\Qq_{645213}\,.
\end{split}
\end{equation}
\begin{equation}
\begin{split}
\Delta \Qq_{{\begin{picture}(6,5)
\put(1,3){\circle*{0.7}}
\put(2,2){\circle*{0.7}}
\put(1,3){\Line(1,-1)}
\put(3,4){\circle*{0.7}}
\put(4,2){\circle*{0.7}}
\put(5,3){\circle*{0.7}}
\put(6,2){\circle*{0.7}}
\put(5,3){\Line(-1,-1)}
\put(5,3){\Line(1,-1)}
\put(3,4){\Line(-2,-1)}
\put(3,4){\Line(2,-1)}
\put(3,4){\circle*{1}}
\end{picture}}} \ =\ & \Qq_{{\begin{picture}(6,5)
\put(1,3){\circle*{0.7}}
\put(2,2){\circle*{0.7}}
\put(1,3){\Line(1,-1)}
\put(3,4){\circle*{0.7}}
\put(4,2){\circle*{0.7}}
\put(5,3){\circle*{0.7}}
\put(6,2){\circle*{0.7}}
\put(5,3){\Line(-1,-1)}
\put(5,3){\Line(1,-1)}
\put(3,4){\Line(-2,-1)}
\put(3,4){\Line(2,-1)}
\put(3,4){\circle*{1}}
\end{picture}}}\,\otimes\, 1 \ +\ \Qq_{{\begin{picture}(5,4)
\put(1,2){\circle*{0.7}}
\put(2,1){\circle*{0.7}}
\put(1,2){\Line(1,-1)}
\put(3,3){\circle*{0.7}}
\put(4,1){\circle*{0.7}}
\put(5,2){\circle*{0.7}}
\put(5,2){\Line(-1,-1)}
\put(3,3){\Line(-2,-1)}
\put(3,3){\Line(2,-1)}
\put(3,3){\circle*{1}}
\end{picture}}}\,\otimes\, \Qq_{{\begin{picture}(1,2)
\put(1,1){\circle*{0.7}}
\put(1,1){\circle*{1}}
\end{picture}}} \ +\ \Qq_{{\begin{picture}(4,4)
\put(1,2){\circle*{0.7}}
\put(2,1){\circle*{0.7}}
\put(1,2){\Line(1,-1)}
\put(3,3){\circle*{0.7}}
\put(4,2){\circle*{0.7}}
\put(3,3){\Line(-2,-1)}
\put(3,3){\Line(1,-1)}
\put(3,3){\circle*{1}}
\end{picture}}}\,\otimes\, \Qq_{{\begin{picture}(2,3)
\put(1,2){\circle*{0.7}}
\put(2,1){\circle*{0.7}}
\put(1,2){\Line(1,-1)}
\put(1,2){\circle*{1}}
\end{picture}}} \\
& +\ \Qq_{{\begin{picture}(3,4)
\put(1,2){\circle*{0.7}}
\put(2,1){\circle*{0.7}}
\put(1,2){\Line(1,-1)}
\put(3,3){\circle*{0.7}}
\put(3,3){\Line(-2,-1)}
\put(3,3){\circle*{1}}
\end{picture}}}\,\otimes\, \Qq_{{\begin{picture}(3,3)
\put(1,1){\circle*{0.7}}
\put(2,2){\circle*{0.7}}
\put(3,1){\circle*{0.7}}
\put(2,2){\Line(-1,-1)}
\put(2,2){\Line(1,-1)}
\put(2,2){\circle*{1}}
\end{picture}}} \ +\ \Qq_{{\begin{picture}(2,3)
\put(1,2){\circle*{0.7}}
\put(2,1){\circle*{0.7}}
\put(1,2){\Line(1,-1)}
\put(1,2){\circle*{1}}
\end{picture}}}\,\otimes\, \Qq_{{\begin{picture}(4,5)
\put(1,4){\circle*{0.7}}
\put(2,2){\circle*{0.7}}
\put(3,3){\circle*{0.7}}
\put(4,2){\circle*{0.7}}
\put(3,3){\Line(-1,-1)}
\put(3,3){\Line(1,-1)}
\put(1,4){\Line(2,-1)}
\put(1,4){\circle*{1}}
\end{picture}}} \ +\ \Qq_{{\begin{picture}(1,2)
\put(1,1){\circle*{0.7}}
\put(1,1){\circle*{1}}
\end{picture}}}\,\otimes\, \Qq_{{\begin{picture}(5,5)
\put(1,3){\circle*{0.7}}
\put(2,4){\circle*{0.7}}
\put(3,2){\circle*{0.7}}
\put(4,3){\circle*{0.7}}
\put(5,2){\circle*{0.7}}
\put(4,3){\Line(-1,-1)}
\put(4,3){\Line(1,-1)}
\put(2,4){\Line(-1,-1)}
\put(2,4){\Line(2,-1)}
\put(2,4){\circle*{1}}
\end{picture}}} \\
& +\ 1\,\otimes\, \Qq_{{\begin{picture}(6,5)
\put(1,3){\circle*{0.7}}
\put(2,2){\circle*{0.7}}
\put(1,3){\Line(1,-1)}
\put(3,4){\circle*{0.7}}
\put(4,2){\circle*{0.7}}
\put(5,3){\circle*{0.7}}
\put(6,2){\circle*{0.7}}
\put(5,3){\Line(-1,-1)}
\put(5,3){\Line(1,-1)}
\put(3,4){\Line(-2,-1)}
\put(3,4){\Line(2,-1)}
\put(3,4){\circle*{1}}
\end{picture}}}
\end{split}
\end{equation}

\subsection{Antipode of $\PBT$}

We have already seen that since $\PBT$ is a connected graded bialgebra of
finite dimension in each degree, one can define the antipode $\nu$ of
$\PBT$ without ambiguity. It then endows $\PBT$ with a structure of graded
Hopf algebra.
Its formula on the basis $\Qq_T$ is:

\begin{equation}
\label{antipQ}
\nu(\Qq_T) = \sum_{I\vDash n} {(-1)^k \Qq_{w_T(I,0)} \cdots
\Qq_{w_T(I,k-1)}}\,,
\end{equation}
where $k$ is the length of $I$ and $w(I,j)$ is the restriction of the word $w$
to the alphabet interval $[i_1+\cdots+i_j+1,i_1+\cdots+i_{j+1}]$.

\begin{note}{\rm
In Section~\ref{multiplicative}, we define analogs of the complete homogeneous
symmetric functions denoted by $\HH_T$ and elementary symmetric functions
denoted by $\EE_T$ in $\PBT$.
Then the image of $\HH_T$ by the antipode, where $T$ is a right comb tree, is
equal up to sign to $T$ reversed, on the $\EE$ basis.

One proves this property by remarking that the $\HH_T$ correspond to products
of complete functions of $\NCSF$ whereas the $\EE_T$ corresponding to left
comb trees are products of elementary functions of $\NCSF$.
}
\end{note}

Notice that the antipode is not an involution as one can see:
\begin{equation}
\nu(\nu(\HH_{{\begin{picture}(3,4)\put(1,2){\circle*{0.7}}\put(2,1){\circle*{0.7}}\put(1,2){\Line(1,-1)}\put(3,3){\circle*{0.7}}\put(3,3){\Line(-2,-1)}\put(3,3){\circle*{1}}\end{picture}}}
\ \,  )) = 2
\HH_{{\begin{picture}(3,3)\put(1,1){\circle*{0.7}}\put(2,2){\circle*{0.7}}\put(3,1){\circle*{0.7}}\put(2,2){\Line(-1,-1)}\put(2,2){\Line(1,-1)}\put(2,2){\circle*{1}}\end{picture}}}
-         2
          \HH_{{\begin{picture}(3,4)\put(1,3){\circle*{0.7}}\put(2,1){\circle*{0.7}}\put(3,2){\circle*{0.7}}\put(3,2){\Line(-1,-1)}\put(1,3){\Line(2,-1)}\put(1,3){\circle*{1}}\end{picture}}}
+ \HH_{{\begin{picture}(3,4)\put(1,2){\circle*{0.7}}\put(2,1){\circle*{0.7}}\put(1,2){\Line(1,-1)}\put(3,3){\circle*{0.7}}\put(3,3){\Line(-2,-1)}\put(3,3){\circle*{1}}\end{picture}}}
\ \ .
\end{equation}

\subsection{Pairs of Fomin graphs}

One can build a pair of graded graphs $(\Gamma,\Gamma^*)$ in duality
as in Fomin's setting~\cite{Fo}, whose vertices of degree $n$ are the binary
trees of size $n$. In $\Gamma$, there is an edge between $T$ and $T'$ if
$\Pp_{T'}$ appears in the product $\Pp_T\Pp_\bullet$ (the dot $\bullet$ is
the tree of size $1$). In $\Gamma^*$, this edge appears if $\Qq_{T'}$
appears in $\Qq_T\Qq_\bullet$. The sylvester correspondence is the Fomin
correspondence associated to this pair of graphs.

\begin{figure}[ht]
\begingroup
\newdimen\vcadre\vcadre=0.1cm 
\newdimen\hcadre\hcadre=0.1cm 
\setlength\unitlength{1.2mm}

$\xymatrix@R=1cm@C=3mm{
 &  &  &  &  &  &  & *{\GrTeXBox{\emptyset}}\arx1[d]& \\
 &  &  &  &  &  &  & *{\GrTeXBox{{\begin{picture}(1,2)\put(1,1){\circle*{0.7}}\put(1,1){\circle*{1}}\end{picture}}}}\arx1[rrrd]\arx1[llld]& \\
 &  &  &  & *{\GrTeXBox{{\begin{picture}(2,3)\put(1,1){\circle*{0.7}}\put(2,2){\circle*{0.7}}\put(2,2){\Line(-1,-1)}\put(2,2){\circle*{1}}\end{picture}}}}\arx1[llld]\arx1[ld]&  &  &  &  &  & *{\GrTeXBox{{\begin{picture}(2,3)\put(1,2){\circle*{0.7}}\put(2,1){\circle*{0.7}}\put(1,2){\Line(1,-1)}\put(1,2){\circle*{1}}\end{picture}}}}\arx1[lllld]\arx1[lld]\arx1[rrd]& \\
 & *{\GrTeXBox{{\begin{picture}(3,4)\put(1,1){\circle*{0.7}}\put(2,2){\circle*{0.7}}\put(2,2){\Line(-1,-1)}\put(3,3){\circle*{0.7}}\put(3,3){\Line(-1,-1)}\put(3,3){\circle*{1}}\end{picture}}}}\arx1[d]\arx1[ld]&  & *{\GrTeXBox{{\begin{picture}(3,3)\put(1,1){\circle*{0.7}}\put(2,2){\circle*{0.7}}\put(3,1){\circle*{0.7}}\put(2,2){\Line(-1,-1)}\put(2,2){\Line(1,-1)}\put(2,2){\circle*{1}}\end{picture}}}}\arx1[rd]\arx1[d]\arx1[ld]&  &  & *{\GrTeXBox{{\begin{picture}(3,4)\put(1,2){\circle*{0.7}}\put(2,1){\circle*{0.7}}\put(1,2){\Line(1,-1)}\put(3,3){\circle*{0.7}}\put(3,3){\Line(-2,-1)}\put(3,3){\circle*{1}}\end{picture}}}}\arx1[d]\arx1[ld]&  & *{\GrTeXBox{{\begin{picture}(3,4)\put(1,3){\circle*{0.7}}\put(2,1){\circle*{0.7}}\put(3,2){\circle*{0.7}}\put(3,2){\Line(-1,-1)}\put(1,3){\Line(2,-1)}\put(1,3){\circle*{1}}\end{picture}}}}\arx1[rd]\arx1[d]\arx1[ld]&  &  &  & *{\GrTeXBox{{\begin{picture}(3,4)\put(1,3){\circle*{0.7}}\put(2,2){\circle*{0.7}}\put(3,1){\circle*{0.7}}\put(2,2){\Line(1,-1)}\put(1,3){\Line(1,-1)}\put(1,3){\circle*{1}}\end{picture}}}}\arx1[rd]\arx1[d]\arx1[ld]\arx1[lld]& \\
*{\GrTeXBox{{\begin{picture}(4,5)\put(1,1){\circle*{0.7}}\put(2,2){\circle*{0.7}}\put(2,2){\Line(-1,-1)}\put(3,3){\circle*{0.7}}\put(3,3){\Line(-1,-1)}\put(4,4){\circle*{0.7}}\put(4,4){\Line(-1,-1)}\put(4,4){\circle*{1}}\end{picture}}}}& *{\GrTeXBox{{\begin{picture}(4,4)\put(1,1){\circle*{0.7}}\put(2,2){\circle*{0.7}}\put(2,2){\Line(-1,-1)}\put(3,3){\circle*{0.7}}\put(4,2){\circle*{0.7}}\put(3,3){\Line(-1,-1)}\put(3,3){\Line(1,-1)}\put(3,3){\circle*{1}}\end{picture}}}}& *{\GrTeXBox{{\begin{picture}(4,5)\put(1,2){\circle*{0.7}}\put(2,3){\circle*{0.7}}\put(3,2){\circle*{0.7}}\put(2,3){\Line(-1,-1)}\put(2,3){\Line(1,-1)}\put(4,4){\circle*{0.7}}\put(4,4){\Line(-2,-1)}\put(4,4){\circle*{1}}\end{picture}}}}& *{\GrTeXBox{{\begin{picture}(4,4)\put(1,2){\circle*{0.7}}\put(2,3){\circle*{0.7}}\put(3,1){\circle*{0.7}}\put(4,2){\circle*{0.7}}\put(4,2){\Line(-1,-1)}\put(2,3){\Line(-1,-1)}\put(2,3){\Line(2,-1)}\put(2,3){\circle*{1}}\end{picture}}}}& *{\GrTeXBox{{\begin{picture}(4,4)\put(1,2){\circle*{0.7}}\put(2,3){\circle*{0.7}}\put(3,2){\circle*{0.7}}\put(4,1){\circle*{0.7}}\put(3,2){\Line(1,-1)}\put(2,3){\Line(-1,-1)}\put(2,3){\Line(1,-1)}\put(2,3){\circle*{1}}\end{picture}}}}& *{\GrTeXBox{{\begin{picture}(4,5)\put(1,2){\circle*{0.7}}\put(2,1){\circle*{0.7}}\put(1,2){\Line(1,-1)}\put(3,3){\circle*{0.7}}\put(3,3){\Line(-2,-1)}\put(4,4){\circle*{0.7}}\put(4,4){\Line(-1,-1)}\put(4,4){\circle*{1}}\end{picture}}}}& *{\GrTeXBox{{\begin{picture}(4,4)\put(1,2){\circle*{0.7}}\put(2,1){\circle*{0.7}}\put(1,2){\Line(1,-1)}\put(3,3){\circle*{0.7}}\put(4,2){\circle*{0.7}}\put(3,3){\Line(-2,-1)}\put(3,3){\Line(1,-1)}\put(3,3){\circle*{1}}\end{picture}}}}& *{\GrTeXBox{{\begin{picture}(4,5)\put(1,3){\circle*{0.7}}\put(2,1){\circle*{0.7}}\put(3,2){\circle*{0.7}}\put(3,2){\Line(-1,-1)}\put(1,3){\Line(2,-1)}\put(4,4){\circle*{0.7}}\put(4,4){\Line(-3,-1)}\put(4,4){\circle*{1}}\end{picture}}}}& *{\GrTeXBox{{\begin{picture}(4,5)\put(1,4){\circle*{0.7}}\put(2,1){\circle*{0.7}}\put(3,2){\circle*{0.7}}\put(3,2){\Line(-1,-1)}\put(4,3){\circle*{0.7}}\put(4,3){\Line(-1,-1)}\put(1,4){\Line(3,-1)}\put(1,4){\circle*{1}}\end{picture}}}}& *{\GrTeXBox{{\begin{picture}(4,5)\put(1,4){\circle*{0.7}}\put(2,2){\circle*{0.7}}\put(3,3){\circle*{0.7}}\put(4,2){\circle*{0.7}}\put(3,3){\Line(-1,-1)}\put(3,3){\Line(1,-1)}\put(1,4){\Line(2,-1)}\put(1,4){\circle*{1}}\end{picture}}}}& *{\GrTeXBox{{\begin{picture}(4,5)\put(1,3){\circle*{0.7}}\put(2,2){\circle*{0.7}}\put(3,1){\circle*{0.7}}\put(2,2){\Line(1,-1)}\put(1,3){\Line(1,-1)}\put(4,4){\circle*{0.7}}\put(4,4){\Line(-3,-1)}\put(4,4){\circle*{1}}\end{picture}}}}& *{\GrTeXBox{{\begin{picture}(4,5)\put(1,4){\circle*{0.7}}\put(2,2){\circle*{0.7}}\put(3,1){\circle*{0.7}}\put(2,2){\Line(1,-1)}\put(4,3){\circle*{0.7}}\put(4,3){\Line(-2,-1)}\put(1,4){\Line(3,-1)}\put(1,4){\circle*{1}}\end{picture}}}}& *{\GrTeXBox{{\begin{picture}(4,5)\put(1,4){\circle*{0.7}}\put(2,3){\circle*{0.7}}\put(3,1){\circle*{0.7}}\put(4,2){\circle*{0.7}}\put(4,2){\Line(-1,-1)}\put(2,3){\Line(2,-1)}\put(1,4){\Line(1,-1)}\put(1,4){\circle*{1}}\end{picture}}}}& *{\GrTeXBox{{\begin{picture}(4,5)\put(1,4){\circle*{0.7}}\put(2,3){\circle*{0.7}}\put(3,2){\circle*{0.7}}\put(4,1){\circle*{0.7}}\put(3,2){\Line(1,-1)}\put(2,3){\Line(1,-1)}\put(1,4){\Line(1,-1)}\put(1,4){\circle*{1}}\end{picture}}}}& \\
}$

\vskip1cm

$\xymatrix@R=1cm@C=3mm{
 &  &  &  &  &  &  & *{\GrTeXBox{\emptyset}}\arx1[d]& \\
 &  &  &  &  &  &  & *{\GrTeXBox{{\begin{picture}(1,2)\put(1,1){\circle*{0.7}}\put(1,1){\circle*{1}}\end{picture}}}}\arx1[rrrd]\arx1[llld]& \\
 &  &  &  & *{\GrTeXBox{{\begin{picture}(2,3)\put(1,1){\circle*{0.7}}\put(2,2){\circle*{0.7}}\put(2,2){\Line(-1,-1)}\put(2,2){\circle*{1}}\end{picture}}}}\arx1[llld]\arx1[rrrrd]\arx1[ld]&  &  &  &  &  & *{\GrTeXBox{{\begin{picture}(2,3)\put(1,2){\circle*{0.7}}\put(2,1){\circle*{0.7}}\put(1,2){\Line(1,-1)}\put(1,2){\circle*{1}}\end{picture}}}}\arx1[lllld]\arx1[rrd]\arx1[llllllld]& \\
 & *{\GrTeXBox{{\begin{picture}(3,4)\put(1,1){\circle*{0.7}}\put(2,2){\circle*{0.7}}\put(2,2){\Line(-1,-1)}\put(3,3){\circle*{0.7}}\put(3,3){\Line(-1,-1)}\put(3,3){\circle*{1}}\end{picture}}}}\arx1[d]\arx1[ld]\arx1[rrrrrrrd]\arx1[rrd]&  & *{\GrTeXBox{{\begin{picture}(3,3)\put(1,1){\circle*{0.7}}\put(2,2){\circle*{0.7}}\put(3,1){\circle*{0.7}}\put(2,2){\Line(-1,-1)}\put(2,2){\Line(1,-1)}\put(2,2){\circle*{1}}\end{picture}}}}\arx1[lld]\arx1[rrrrrrd]\arx1[rd]\arx1[ld]&  &  & *{\GrTeXBox{{\begin{picture}(3,4)\put(1,2){\circle*{0.7}}\put(2,1){\circle*{0.7}}\put(1,2){\Line(1,-1)}\put(3,3){\circle*{0.7}}\put(3,3){\Line(-2,-1)}\put(3,3){\circle*{1}}\end{picture}}}}\arx1[d]\arx1[ld]\arx1[rrrrrd]\arx1[llld]&  & *{\GrTeXBox{{\begin{picture}(3,4)\put(1,3){\circle*{0.7}}\put(2,1){\circle*{0.7}}\put(3,2){\circle*{0.7}}\put(3,2){\Line(-1,-1)}\put(1,3){\Line(2,-1)}\put(1,3){\circle*{1}}\end{picture}}}}\arx1[lld]\arx1[ld]\arx1[rrrrd]\arx1[llllld]&  &  &  & *{\GrTeXBox{{\begin{picture}(3,4)\put(1,3){\circle*{0.7}}\put(2,2){\circle*{0.7}}\put(3,1){\circle*{0.7}}\put(2,2){\Line(1,-1)}\put(1,3){\Line(1,-1)}\put(1,3){\circle*{1}}\end{picture}}}}\arx1[lllllld]\arx1[rd]\arx1[lld]\arx1[lllllllld]& \\
*{\GrTeXBox{{\begin{picture}(4,5)\put(1,1){\circle*{0.7}}\put(2,2){\circle*{0.7}}\put(2,2){\Line(-1,-1)}\put(3,3){\circle*{0.7}}\put(3,3){\Line(-1,-1)}\put(4,4){\circle*{0.7}}\put(4,4){\Line(-1,-1)}\put(4,4){\circle*{1}}\end{picture}}}}& *{\GrTeXBox{{\begin{picture}(4,4)\put(1,1){\circle*{0.7}}\put(2,2){\circle*{0.7}}\put(2,2){\Line(-1,-1)}\put(3,3){\circle*{0.7}}\put(4,2){\circle*{0.7}}\put(3,3){\Line(-1,-1)}\put(3,3){\Line(1,-1)}\put(3,3){\circle*{1}}\end{picture}}}}& *{\GrTeXBox{{\begin{picture}(4,5)\put(1,2){\circle*{0.7}}\put(2,3){\circle*{0.7}}\put(3,2){\circle*{0.7}}\put(2,3){\Line(-1,-1)}\put(2,3){\Line(1,-1)}\put(4,4){\circle*{0.7}}\put(4,4){\Line(-2,-1)}\put(4,4){\circle*{1}}\end{picture}}}}& *{\GrTeXBox{{\begin{picture}(4,4)\put(1,2){\circle*{0.7}}\put(2,3){\circle*{0.7}}\put(3,1){\circle*{0.7}}\put(4,2){\circle*{0.7}}\put(4,2){\Line(-1,-1)}\put(2,3){\Line(-1,-1)}\put(2,3){\Line(2,-1)}\put(2,3){\circle*{1}}\end{picture}}}}& *{\GrTeXBox{{\begin{picture}(4,4)\put(1,2){\circle*{0.7}}\put(2,3){\circle*{0.7}}\put(3,2){\circle*{0.7}}\put(4,1){\circle*{0.7}}\put(3,2){\Line(1,-1)}\put(2,3){\Line(-1,-1)}\put(2,3){\Line(1,-1)}\put(2,3){\circle*{1}}\end{picture}}}}& *{\GrTeXBox{{\begin{picture}(4,5)\put(1,2){\circle*{0.7}}\put(2,1){\circle*{0.7}}\put(1,2){\Line(1,-1)}\put(3,3){\circle*{0.7}}\put(3,3){\Line(-2,-1)}\put(4,4){\circle*{0.7}}\put(4,4){\Line(-1,-1)}\put(4,4){\circle*{1}}\end{picture}}}}& *{\GrTeXBox{{\begin{picture}(4,4)\put(1,2){\circle*{0.7}}\put(2,1){\circle*{0.7}}\put(1,2){\Line(1,-1)}\put(3,3){\circle*{0.7}}\put(4,2){\circle*{0.7}}\put(3,3){\Line(-2,-1)}\put(3,3){\Line(1,-1)}\put(3,3){\circle*{1}}\end{picture}}}}& *{\GrTeXBox{{\begin{picture}(4,5)\put(1,3){\circle*{0.7}}\put(2,1){\circle*{0.7}}\put(3,2){\circle*{0.7}}\put(3,2){\Line(-1,-1)}\put(1,3){\Line(2,-1)}\put(4,4){\circle*{0.7}}\put(4,4){\Line(-3,-1)}\put(4,4){\circle*{1}}\end{picture}}}}& *{\GrTeXBox{{\begin{picture}(4,5)\put(1,4){\circle*{0.7}}\put(2,1){\circle*{0.7}}\put(3,2){\circle*{0.7}}\put(3,2){\Line(-1,-1)}\put(4,3){\circle*{0.7}}\put(4,3){\Line(-1,-1)}\put(1,4){\Line(3,-1)}\put(1,4){\circle*{1}}\end{picture}}}}& *{\GrTeXBox{{\begin{picture}(4,5)\put(1,4){\circle*{0.7}}\put(2,2){\circle*{0.7}}\put(3,3){\circle*{0.7}}\put(4,2){\circle*{0.7}}\put(3,3){\Line(-1,-1)}\put(3,3){\Line(1,-1)}\put(1,4){\Line(2,-1)}\put(1,4){\circle*{1}}\end{picture}}}}& *{\GrTeXBox{{\begin{picture}(4,5)\put(1,3){\circle*{0.7}}\put(2,2){\circle*{0.7}}\put(3,1){\circle*{0.7}}\put(2,2){\Line(1,-1)}\put(1,3){\Line(1,-1)}\put(4,4){\circle*{0.7}}\put(4,4){\Line(-3,-1)}\put(4,4){\circle*{1}}\end{picture}}}}& *{\GrTeXBox{{\begin{picture}(4,5)\put(1,4){\circle*{0.7}}\put(2,2){\circle*{0.7}}\put(3,1){\circle*{0.7}}\put(2,2){\Line(1,-1)}\put(4,3){\circle*{0.7}}\put(4,3){\Line(-2,-1)}\put(1,4){\Line(3,-1)}\put(1,4){\circle*{1}}\end{picture}}}}& *{\GrTeXBox{{\begin{picture}(4,5)\put(1,4){\circle*{0.7}}\put(2,3){\circle*{0.7}}\put(3,1){\circle*{0.7}}\put(4,2){\circle*{0.7}}\put(4,2){\Line(-1,-1)}\put(2,3){\Line(2,-1)}\put(1,4){\Line(1,-1)}\put(1,4){\circle*{1}}\end{picture}}}}& *{\GrTeXBox{{\begin{picture}(4,5)\put(1,4){\circle*{0.7}}\put(2,3){\circle*{0.7}}\put(3,2){\circle*{0.7}}\put(4,1){\circle*{0.7}}\put(3,2){\Line(1,-1)}\put(2,3){\Line(1,-1)}\put(1,4){\Line(1,-1)}\put(1,4){\circle*{1}}\end{picture}}}}& \\
}$
\bigskip

\endgroup

\caption{The two graded graphs in duality.}
\end{figure}
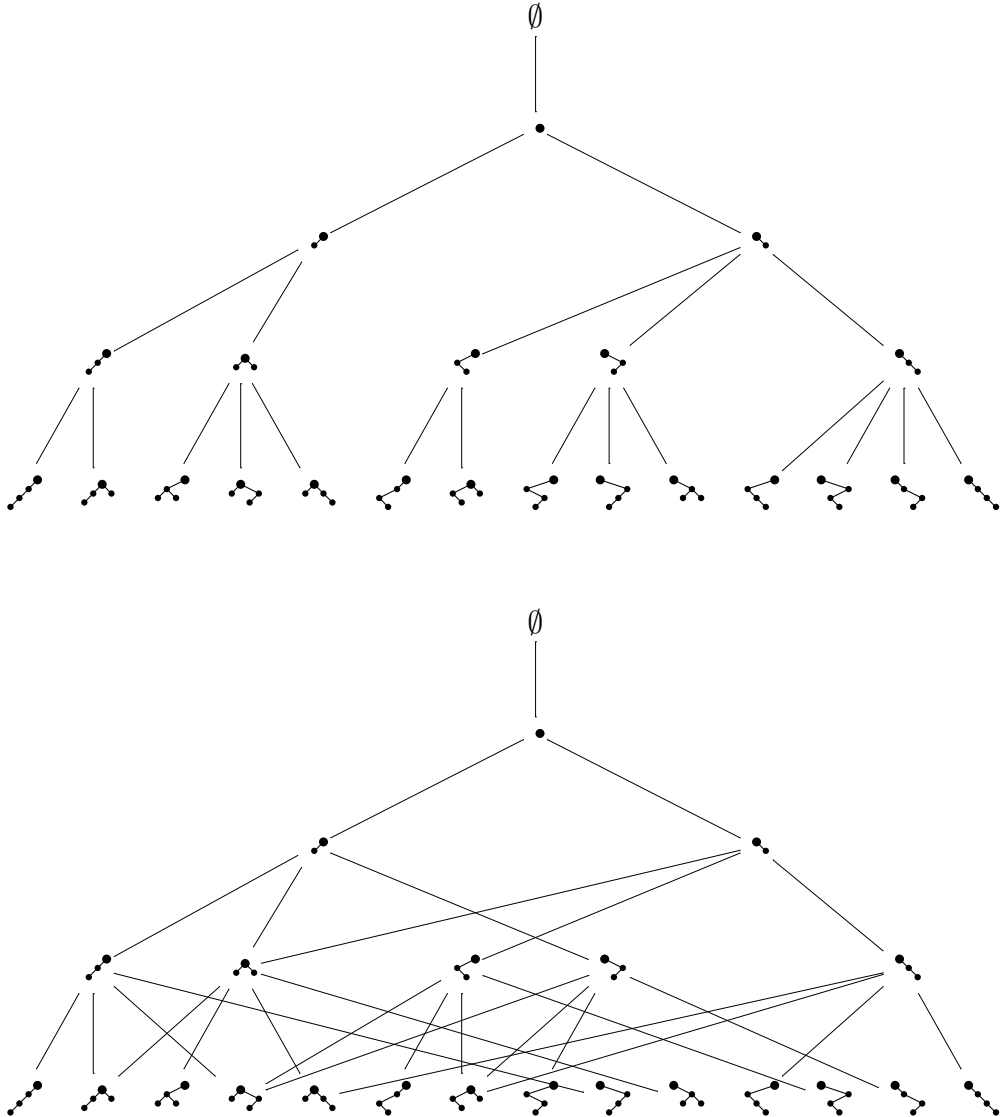

\subsection{The Tamari order and equivalent orders}
\label{Tamari}

In~\cite{LR2}, Loday and Ronco describe the product $\Pp_{T'}\Pp_{T''}$ as
an interval of the so-called \emph{Tamari} order.
The situation is the following: given the structure of all sylvester classes
inside the permutohedron, we can easily prove in our setting that the product
of two $\Pp_{T}$ functions is an interval of the permutohedron (see
Note~\ref{intercombi}). It happens that the restriction of the weak
order to sylvester classes is the same as the Tamari order (see
Theorem~\ref{thmsylvTam}), which yields in particular a  simple  proof of
the result of~\cite{LR2}.

Let us first give some definitions.
Following Stanley in~\cite{Stan2} (ex. 6.32.a p.~234), we define the Tamari
order $O_n$ as the poset of all integer vectors $(a_1,\ldots,a_n)$ such that
$i\leq a_i\leq n$ and such that, if $i\leq j\leq a_i$ then $a_j\leq a_i$,
ordered coordinatewise (see Figures~\ref{tam3} and~\ref{tam4}).

\smallskip
\ifdessin
\begin{figure}[ht]
\newdimen\vcadre\vcadre=0.1cm 
\newdimen\hcadre\hcadre=0.1cm 
$\xymatrix@R=0.5cm@C=8mm{
 & *{\GrTeXBox{\begin{picture}(3,4)\put(1,1){\circle*{0.7}}\put(2,2){\circle*{0.7}}\put(2,2){\Line(-1,-1)}\put(3,3){\circle*{0.7}}\put(3,3){\Line(-1,-1)}\put(3,3){\circle*{1}}\end{picture}}}\arx1[ld]\arx1[rdd] \\
*{\GrTeXBox{\begin{picture}(3,4)\put(1,2){\circle*{0.7}}\put(2,1){\circle*{0.7}}\put(1,2){\Line(1,-1)}\put(3,3){\circle*{0.7}}\put(3,3){\Line(-2,-1)}\put(3,3){\circle*{1}}\end{picture}}}\arx1[d] \\
*{\GrTeXBox{\begin{picture}(3,4)\put(1,3){\circle*{0.7}}\put(2,1){\circle*{0.7}}\put(3,2){\circle*{0.7}}\put(3,2){\Line(-1,-1)}\put(1,3){\Line(2,-1)}\put(1,3){\circle*{1}}\end{picture}}}\arx1[rd]&  & *{\GrTeXBox{\begin{picture}(3,3)\put(1,1){\circle*{0.7}}\put(2,2){\circle*{0.7}}\put(3,1){\circle*{0.7}}\put(2,2){\Line(-1,-1)}\put(2,2){\Line(1,-1)}\put(2,2){\circle*{1}}\end{picture}}}\arx1[ld] \\
 & *{\GrTeXBox{\begin{picture}(3,4)\put(1,3){\circle*{0.7}}\put(2,2){\circle*{0.7}}\put(3,1){\circle*{0.7}}\put(2,2){\Line(1,-1)}\put(1,3){\Line(1,-1)}\put(1,3){\circle*{1}}\end{picture}}} \\
}$
\hfill
\newdimen\vcadre\vcadre=0.2cm 
\newdimen\hcadre\hcadre=0.2cm 
$\xymatrix@R=0.9cm@C=7mm{
 & *{\GrTeXBox{123}}\arx1[ld]\arx1[rdd]\\
*{\GrTeXBox{213}}\arx1[d]\\
*{\GrTeXBox{231}}\arx1[rd]&  & *{\GrTeXBox{312}}\arx1[ld]\\
 & *{\GrTeXBox{321}}\\
}$
\hfill
\newdimen\vcadre\vcadre=0.2cm 
\newdimen\hcadre\hcadre=0.2cm 
$\xymatrix@R=0.9cm@C=7mm{
 & *{\GrTeXBox{111}}\arx1[ld]\arx1[rdd] \\
*{\GrTeXBox{121}}\arx1[d] \\
*{\GrTeXBox{122}}\arx1[rd]&  & *{\GrTeXBox{113}}\arx1[ld] \\
 & *{\GrTeXBox{123}} \\
}$
\caption{\label{tam3}The same order on trees, canonical words, and
Tamari elements of size $3$.}
\end{figure}
\fi

\ifdessin
\begin{figure}[ht]
{\centerline{
\tiny
\rotateleft{
\newdimen\vcadre\vcadre=0.15cm 
\newdimen\hcadre\hcadre=0.15cm 
\setlength\unitlength{1.5mm}
$\xymatrix@R=0.5cm@C=4mm{
 &  & *{\GrTeXBox{\begin{picture}(4,5)\put(1,1){\circle*{0.7}}\put(2,2){\circle*{0.7}}\put(2,2){\Line(-1,-1)}\put(3,3){\circle*{0.7}}\put(3,3){\Line(-1,-1)}\put(4,4){\circle*{0.7}}\put(4,4){\Line(-1,-1)}\put(4,4){\circle*{1}}\end{picture}}}\arx1[ld]\arx1[dd]\arx1[rrddd]& \\
 & *{\GrTeXBox{\begin{picture}(4,5)\put(1,2){\circle*{0.7}}\put(2,1){\circle*{0.7}}\put(1,2){\Line(1,-1)}\put(3,3){\circle*{0.7}}\put(3,3){\Line(-2,-1)}\put(4,4){\circle*{0.7}}\put(4,4){\Line(-1,-1)}\put(4,4){\circle*{1}}\end{picture}}}\arx1[ld]\arx1[rddd]& \\
*{\GrTeXBox{\begin{picture}(4,5)\put(1,3){\circle*{0.7}}\put(2,1){\circle*{0.7}}\put(3,2){\circle*{0.7}}\put(3,2){\Line(-1,-1)}\put(1,3){\Line(2,-1)}\put(4,4){\circle*{0.7}}\put(4,4){\Line(-3,-1)}\put(4,4){\circle*{1}}\end{picture}}}\arx1[d]\arx1[rd]&  & *{\GrTeXBox{\begin{picture}(4,5)\put(1,2){\circle*{0.7}}\put(2,3){\circle*{0.7}}\put(3,2){\circle*{0.7}}\put(2,3){\Line(-1,-1)}\put(2,3){\Line(1,-1)}\put(4,4){\circle*{0.7}}\put(4,4){\Line(-2,-1)}\put(4,4){\circle*{1}}\end{picture}}}\arx1[ld]\arx1[rdd]& \\
*{\GrTeXBox{\begin{picture}(4,5)\put(1,4){\circle*{0.7}}\put(2,1){\circle*{0.7}}\put(3,2){\circle*{0.7}}\put(3,2){\Line(-1,-1)}\put(4,3){\circle*{0.7}}\put(4,3){\Line(-1,-1)}\put(1,4){\Line(3,-1)}\put(1,4){\circle*{1}}\end{picture}}}\arx1[d]\arx1[rrdd]& *{\GrTeXBox{\begin{picture}(4,5)\put(1,3){\circle*{0.7}}\put(2,2){\circle*{0.7}}\put(3,1){\circle*{0.7}}\put(2,2){\Line(1,-1)}\put(1,3){\Line(1,-1)}\put(4,4){\circle*{0.7}}\put(4,4){\Line(-3,-1)}\put(4,4){\circle*{1}}\end{picture}}}\arx1[ld]&  &  & *{\GrTeXBox{\begin{picture}(4,4)\put(1,1){\circle*{0.7}}\put(2,2){\circle*{0.7}}\put(2,2){\Line(-1,-1)}\put(3,3){\circle*{0.7}}\put(4,2){\circle*{0.7}}\put(3,3){\Line(-1,-1)}\put(3,3){\Line(1,-1)}\put(3,3){\circle*{1}}\end{picture}}}\arx1[lld]\arx1[dd]& \\
*{\GrTeXBox{\begin{picture}(4,5)\put(1,4){\circle*{0.7}}\put(2,2){\circle*{0.7}}\put(3,1){\circle*{0.7}}\put(2,2){\Line(1,-1)}\put(4,3){\circle*{0.7}}\put(4,3){\Line(-2,-1)}\put(1,4){\Line(3,-1)}\put(1,4){\circle*{1}}\end{picture}}}\arx1[rd]&  & *{\GrTeXBox{\begin{picture}(4,4)\put(1,2){\circle*{0.7}}\put(2,1){\circle*{0.7}}\put(1,2){\Line(1,-1)}\put(3,3){\circle*{0.7}}\put(4,2){\circle*{0.7}}\put(3,3){\Line(-2,-1)}\put(3,3){\Line(1,-1)}\put(3,3){\circle*{1}}\end{picture}}}\arx1[d]& *{\GrTeXBox{\begin{picture}(4,4)\put(1,2){\circle*{0.7}}\put(2,3){\circle*{0.7}}\put(3,1){\circle*{0.7}}\put(4,2){\circle*{0.7}}\put(4,2){\Line(-1,-1)}\put(2,3){\Line(-1,-1)}\put(2,3){\Line(2,-1)}\put(2,3){\circle*{1}}\end{picture}}}\arx1[lld]\arx1[rd]& \\
 & *{\GrTeXBox{\begin{picture}(4,5)\put(1,4){\circle*{0.7}}\put(2,3){\circle*{0.7}}\put(3,1){\circle*{0.7}}\put(4,2){\circle*{0.7}}\put(4,2){\Line(-1,-1)}\put(2,3){\Line(2,-1)}\put(1,4){\Line(1,-1)}\put(1,4){\circle*{1}}\end{picture}}}\arx1[rd]& *{\GrTeXBox{\begin{picture}(4,5)\put(1,4){\circle*{0.7}}\put(2,2){\circle*{0.7}}\put(3,3){\circle*{0.7}}\put(4,2){\circle*{0.7}}\put(3,3){\Line(-1,-1)}\put(3,3){\Line(1,-1)}\put(1,4){\Line(2,-1)}\put(1,4){\circle*{1}}\end{picture}}}\arx1[d]&  & *{\GrTeXBox{\begin{picture}(4,4)\put(1,2){\circle*{0.7}}\put(2,3){\circle*{0.7}}\put(3,2){\circle*{0.7}}\put(4,1){\circle*{0.7}}\put(3,2){\Line(1,-1)}\put(2,3){\Line(-1,-1)}\put(2,3){\Line(1,-1)}\put(2,3){\circle*{1}}\end{picture}}}\arx1[lld]& \\
 &  & *{\GrTeXBox{\begin{picture}(4,5)\put(1,4){\circle*{0.7}}\put(2,3){\circle*{0.7}}\put(3,2){\circle*{0.7}}\put(4,1){\circle*{0.7}}\put(3,2){\Line(1,-1)}\put(2,3){\Line(1,-1)}\put(1,4){\Line(1,-1)}\put(1,4){\circle*{1}}\end{picture}}}& \\
}$
\hfill
\newdimen\vcadre\vcadre=0.2cm 
\newdimen\hcadre\hcadre=0.2cm 
$\xymatrix@R=1.0cm@C=4mm{
 &  & *{\GrTeXBox{1234}}\arx1[ld]\arx1[dd]\arx1[rrddd]& \\
 & *{\GrTeXBox{2134}}\arx1[ld]\arx1[rddd]& \\
*{\GrTeXBox{2314}}\arx1[d]\arx1[rd]&  & *{\GrTeXBox{3124}}\arx1[ld]\arx1[rdd]& \\
*{\GrTeXBox{2341}}\arx1[d]\arx1[rrdd]& *{\GrTeXBox{3214}}\arx1[ld]&  &  & *{\GrTeXBox{4123}}\arx1[lld]\arx1[dd]& \\
*{\GrTeXBox{3241}}\arx1[rd]&  & *{\GrTeXBox{4213}}\arx1[d]& *{\GrTeXBox{3412}}\arx1[lld]\arx1[rd]& \\
 & *{\GrTeXBox{3421}}\arx1[rd]& *{\GrTeXBox{4231}}\arx1[d]&  & *{\GrTeXBox{4312}}\arx1[lld]& \\
 &  & *{\GrTeXBox{4321}}& \\
}$
\hfill
\newdimen\vcadre\vcadre=0.2cm 
\newdimen\hcadre\hcadre=0.2cm 
$\xymatrix@R=1.0cm@C=4mm{
 &  & *{\GrTeXBox{1111}}\arx1[ld]\arx1[dd]\arx1[rrddd]& \\
 & *{\GrTeXBox{1211}}\arx1[ld]\arx1[rddd]& \\
*{\GrTeXBox{1221}}\arx1[d]\arx1[rd]&  & *{\GrTeXBox{1131}}\arx1[ld]\arx1[rdd]& \\
*{\GrTeXBox{1222}}\arx1[d]\arx1[rrdd]& *{\GrTeXBox{1231}}\arx1[ld]&  &  & *{\GrTeXBox{1114}}\arx1[lld]\arx1[dd]& \\
*{\GrTeXBox{1232}}\arx1[rd]&  & *{\GrTeXBox{1214}}\arx1[d]& *{\GrTeXBox{1133}}\arx1[lld]\arx1[rd]& \\
 & *{\GrTeXBox{1233}}\arx1[rd]& *{\GrTeXBox{1224}}\arx1[d]&  & *{\GrTeXBox{1134}}\arx1[lld]& \\
 &  & *{\GrTeXBox{1234}}& \\
}$
}
}}
\caption{\label{tam4}The same order on trees, canonical words, and Tamari
elements of size $4$.}
\end{figure}
\fi

Let $O'_n$ be the \emph{sylvestrohedron order} defined on sylvester
classes as follows: a sylvester class $S$ is smaller than $S'$ if there exist
$\sigma\in S$ and $\sigma'\in S'$ such that $\sigma<\sigma'$ for the
weak order.

\begin{thm}
\label{thmsylvTam}
The sylvestrohedron order coincides with the Tamari order.
\end{thm}

To prove this property, we will go through a third equivalent order. Let
$O''_n$ be the \emph{sylvester order} defined on canonical sylvester words as
the restriction of the weak order to those elements.

\begin{pf}
The proof of Theorem~\ref{thmsylvTam} results from the following three lemmas:

\begin{lem}
\label{comparsylv}
Let $\sigma$ and $\sigma'$ be two permutations such that $\sigma$ is smaller
than $\sigma'$ for the weak order.
Then the canonical word corresponding to $\sigma$ is smaller than or equal to
the canonical word of $\sigma'$.
\end{lem}

\begin{lem}
\label{sylv2canon}
The sylvestrohedron order coincides with the sylvester order.
\end{lem}

\begin{lem}
\label{tamari2canon}
The Tamari order coincides with the sylvester order.
\end{lem}

Let us first prove Lemma~\ref{comparsylv}. We only need to prove it for an
elementary transposition.
Assume that an elementary transposition $s_i$ sends a permutation $\sigma$ to
$\sigma'$, belonging to another sylvester class. Let us prove that there is an
elementary transposition that sends the canonical word associated with
$\sigma$ inside the class of $\sigma'$.

If $\sigma$ is the canonical word of its class, then it is smaller than the
canonical word of the class of $\sigma'$ by transitivity. Assume now
that $\sigma$ is not a canonical word.
Then there is another elementary transposition $s_j$ that sends $\sigma$ to
$\sigma_0$, belonging to the same sylvester class (let us recall that thanks
to~\cite{BW2}, the sylvester classes are intervals of the permutohedron). If
$j\not=i-1$ and $j\not=i+1$, then $s_i$ sends $\sigma_0$ to an element of the
same sylvester class as $\sigma'$. Otherwise, let us assume that $j=i-1$.
Apply $s_i$ and then $s_{i-1}$ to $\sigma_0$. The resulting element is in the
same sylvester class as $\sigma'$. Applying this property until $\sigma_0$ is
a canonical word proves that there exists an elementary transposition that
sends the canonical word associated with $\sigma$ to the class of $\sigma'$.
\qed

\smallskip
Let us now prove Lemma~\ref{sylv2canon}. The isomorphism between
both orders is trivial: a sylvester class is sent to its canonical word.
Now, by definition, if $\sigma<\sigma'$ for the sylvester order, then
$\sigma<\sigma'$ for the sylvestrohedron order.
Lemma~\ref{comparsylv} proves the converse.
\qed

Let us finally prove Lemma~\ref{tamari2canon}.
Notice that there exists a well-known simple bijection between the elements of
the Tamari poset as defined before and the canonical sylvester words.
Indeed, send each permutation to the sequence defined as follows: for each $i$,
compute the number of elements smaller than $i$ and to the right of $i$. Then
add 1 to each component of the resulting vector. It is a Tamari element.
Conversely, subtract 1 to each component of a Tamari element and rebuild the
permutation which has these numbers of inversions.
For example, if $\sigma=435216$, one finds the sequence $012320$ and the
Tamari element $t=123431$.

Now, given the bijection, it is immediate to see that if two permutations are
comparable for the weak order, then so are the corresponding Tamari
elements for the Tamari order. And conversely, if two Tamari elements are
comparable for the Tamari order, then so are the corresponding permutations.
\qed

Lemmas~\ref{sylv2canon} and~\ref{tamari2canon} together imply the theorem.
\qed
\end{pf}

\begin{note}{\rm
\label{FSYmbug}
The construction of the sylvestrohedron order and the fact that it is
the same as the restriction of the weak order to canonical words does not work
for the other known interesting example, that is the algebra of Free
Symmetric Functions $\FSym$: in this algebra, if one says that a plactic class
is smaller than another one if there is an element of the first one smaller
than an element of the second one for the weak order, this relation is not
transitive.
}
\end{note}

\begin{note}
\label{intercombi}{\rm
It is well known that the set of permutations arising in the shifted shuffle
of two permutohedron intervals is a permutohedron interval.
Since each sylvester class is an interval of the permutohedron and since the
product of quasi-ribbon functions is given by the shifted shuffle, it
immediately comes that the product $\Pp_{T'}\Pp_{T''}$ expressed on the
$\F_\sigma$'s is an interval of the permutohedron, and so, is an interval of
the sylvestrohedron. This property holds for other quotients of
$\FQSym^*$ as soon as the congruence is compatible with de-standardization.
}
\end{note}

\begin{note}{\rm
\label{couvcanons}
As already proved, the Tamari order is the same as the restriction of the
weak order to canonical words. One can translate on canonical words the
covering relation built by Loday and Ronco: given a canonical word $\sigma$,
for any rise $\sigma(i)<\sigma(i+1)$, one builds the permutation obtained from
$\sigma$ by exchanging $\sigma(i+1)$ with the element to its left as long as
it is smaller than or equal to $\sigma(i)$. This last permutation is a
canonical word.

For example, if one chooses $\sigma=(12,10,8,9,6,7,4,2,1,3,5,11)$ and $i=11$,
one obtains $(12,10,8,9,6,7,11,4,2,1,3,5)$, which is canonical.
Doing this on all consecutive elements such that $\sigma(i)<\sigma(i+1)$, one
recovers all covering relations of the Tamari order.
}
\end{note}

Lemma~\ref{comparsylv} also proves that

\begin{cor}
\label{corintercombi}
The intervals of the permutohedron starting at the identity permutation and
finishing at a canonical sylvester permutation are unions of sylvester
classes.
\end{cor}

This property will be the main tool for constructing multiplicative bases.

\subsection{Multiplicative bases}
\label{multiplicative}

In their paper~\cite{LR1}, Loday and Ronco build a multiplicative basis of
$\PBT$ by associating with a tree, a function obtained by multiplying the
$\Pp_T$'s obtained by cutting the right subtrees connected by right edges to
the root of $T$.

However, in our setting, there is a more natural and general way to
build multiplicative bases. Let us first fix the notation:
a tree $T$ is said to be smaller than a tree $T'$, and we write $T<T'$ if
$w_T<w_{T'}$ for the sylvester order.

\begin{def}
\label{HEbase}
Let $T$ be a tree. The \emph{complete} ($\HH$) and \emph{elementary} ($\EE$)
functions of $\PBT$ are respectively defined by
\begin{equation}
\HH_T := \sum_{T'\leq T} \Pp_{T'}\,,
\end{equation}
\begin{equation}
\EE_T := \sum_{T'\geq T} \Pp_{T'}\,.
\end{equation}
\end{def}

The names complete and elementary functions have been chosen on purpose: as we
will see later (see~Section~\ref{morphismes}), these are analogs of the
homogeneous complete and elementary symmetric functions.

For example,
\begin{equation}
\HH_{213} = \Pp_{123} + \Pp_{213}: \qquad
\HH_{{\begin{picture}(3,4)
\put(1,2){\circle*{0.7}}
\put(2,1){\circle*{0.7}}
\put(1,2){\Line(1,-1)}
\put(3,3){\circle*{0.7}}
\put(3,3){\Line(-2,-1)}
\put(3,3){\circle*{1}}
\end{picture}}} \ =\ \Pp_{{\begin{picture}(3,4)
\put(1,1){\circle*{0.7}}
\put(2,2){\circle*{0.7}}
\put(2,2){\Line(-1,-1)}
\put(3,3){\circle*{0.7}}
\put(3,3){\Line(-1,-1)}
\put(3,3){\circle*{1}}
\end{picture}}} \ +\ \Pp_{{\begin{picture}(3,4)
\put(1,2){\circle*{0.7}}
\put(2,1){\circle*{0.7}}
\put(1,2){\Line(1,-1)}
\put(3,3){\circle*{0.7}}
\put(3,3){\Line(-2,-1)}
\put(3,3){\circle*{1}}
\end{picture}}}
\ \ ,
\end{equation}

\begin{equation}
\EE_{213} := \Pp_{213} + \Pp_{231} + \Pp_{321}: \qquad
\EE_{{\begin{picture}(3,4)
\put(1,2){\circle*{0.7}}
\put(2,1){\circle*{0.7}}
\put(1,2){\Line(1,-1)}
\put(3,3){\circle*{0.7}}
\put(3,3){\Line(-2,-1)}
\put(3,3){\circle*{1}}
\end{picture}}} \ =\ \Pp_{{\begin{picture}(3,4)
\put(1,2){\circle*{0.7}}
\put(2,1){\circle*{0.7}}
\put(1,2){\Line(1,-1)}
\put(3,3){\circle*{0.7}}
\put(3,3){\Line(-2,-1)}
\put(3,3){\circle*{1}}
\end{picture}}}\ +\ \Pp_{{\begin{picture}(3,4)
\put(1,3){\circle*{0.7}}
\put(2,1){\circle*{0.7}}
\put(3,2){\circle*{0.7}}
\put(3,2){\Line(-1,-1)}
\put(1,3){\Line(2,-1)}
\put(1,3){\circle*{1}}
\end{picture}}}\ +\ \Pp_{{\begin{picture}(3,4)
\put(1,3){\circle*{0.7}}
\put(2,2){\circle*{0.7}}
\put(3,1){\circle*{0.7}}
\put(2,2){\Line(1,-1)}
\put(1,3){\Line(1,-1)}
\put(1,3){\circle*{1}}
\end{picture}}}
\ \ .
\end{equation}

More examples are given at the end of the paper.
For the matrices $M_{\HH,\Pp}$, see Figures~\ref{h2p23} and~\ref{h2p4}.
For $M_{\EE,\Pp}$, see Figures~\ref{e2p23} and~\ref{e2p4}.
For $M_{\EE,\HH}$, see Figures~\ref{e2h23} and~\ref{e2h4}.

\begin{thm}
\label{multiH}
The basis of $\HH$ functions is a multiplicative basis of $\PBT$, whose
product is given by
\begin{equation}
\HH_{T'} \HH_{T''} = \HH_{T}\,,
\end{equation}
where the canonical word $w_T$ of $T$ is obtained by concatenating
$w_{T''}[k]$ with $w_{T'}$ where $k$ is the size of $T'$.
\end{thm}

For example,
\begin{equation}
\HH_{312} \HH_{45213} = \HH_{78564312}: \qquad
\HH_{{\begin{picture}(3,3)
\put(1,1){\circle*{0.7}}
\put(2,2){\circle*{0.7}}
\put(3,1){\circle*{0.7}}
\put(2,2){\Line(-1,-1)}
\put(2,2){\Line(1,-1)}
\put(2,2){\circle*{1}}
\end{picture}}}\, \HH_{{\begin{picture}(5,4)
\put(1,2){\circle*{0.7}}
\put(2,1){\circle*{0.7}}
\put(1,2){\Line(1,-1)}
\put(3,3){\circle*{0.7}}
\put(4,1){\circle*{0.7}}
\put(5,2){\circle*{0.7}}
\put(5,2){\Line(-1,-1)}
\put(3,3){\Line(-2,-1)}
\put(3,3){\Line(2,-1)}
\put(3,3){\circle*{1}}
\end{picture}}} = \HH_{{\begin{picture}(8,8)
\put(1,6){\circle*{0.7}}
\put(2,7){\circle*{0.7}}
\put(3,6){\circle*{0.7}}
\put(4,4){\circle*{0.7}}
\put(5,3){\circle*{0.7}}
\put(4,4){\Line(1,-1)}
\put(6,5){\circle*{0.7}}
\put(7,3){\circle*{0.7}}
\put(8,4){\circle*{0.7}}
\put(8,4){\Line(-1,-1)}
\put(6,5){\Line(-2,-1)}
\put(6,5){\Line(2,-1)}
\put(3,6){\Line(3,-1)}
\put(2,7){\Line(-1,-1)}
\put(2,7){\Line(1,-1)}
\put(2,7){\circle*{1}}
\end{picture}}}
\ \ .
\end{equation}
Notice that this operation coincides with the \emph{over} operation defined by
Loday-Ronco in~\cite{LR1}.
It consists in grafting the tree $T'$ on the right of the rightmost
element of $T''$.

\begin{pf}
As already pointed out in Note~\ref{corintercombi}, the product in $\FQSym$
of an initial interval of $\SG_i$ by an initial interval of $\SG_{n-i}$ gives
rise to an initial interval of $\SG_n$. This proves that the $\HH$'s
form a multiplicative basis. Now, the greatest possible element of this
product is the word $w_T$ defined in the theorem.
\qed
\end{pf}

The same theorem holds for the $\EE$'s, and the proof needs one small change:
the smallest possible element of a given product is not a canonical word but
can be easily rewritten as the one described in the next theorem. One can use
the fact that sylvester classes are invariant through reversion of the
alphabet, considering $(A,>)$ instead of $(A,<)$.
Indeed, if one rewrites everything in terms of the smallest element of each
sylvester class instead of the greatest one, the product of $\EE_{T'}$ by
$\EE_{T''}$ is given by $\EE_T$ where $w_T=w_{T'}\cdot w_{T''}[k]$.

\begin{thm}
\label{multiE}
The basis of $\EE$ functions is a multiplicative basis of $\PBT$, whose
product is given by
\begin{equation}
\EE_{T'} \EE_{T''} = \EE_{T}\,,
\end{equation}
where $T$ is obtained by connecting $T'$ on the left of the left-most element
of $T''$.
\end{thm}

For example,
\begin{equation}
\EE_{312} \EE_{45213} = \EE_{78531246}: \qquad
\EE_{{\begin{picture}(3,3)
\put(1,1){\circle*{0.7}}
\put(2,2){\circle*{0.7}}
\put(3,1){\circle*{0.7}}
\put(2,2){\Line(-1,-1)}
\put(2,2){\Line(1,-1)}
\put(2,2){\circle*{1}}
\end{picture}}} \EE_{{\begin{picture}(5,4)
\put(1,2){\circle*{0.7}}
\put(2,1){\circle*{0.7}}
\put(1,2){\Line(1,-1)}
\put(3,3){\circle*{0.7}}
\put(4,1){\circle*{0.7}}
\put(5,2){\circle*{0.7}}
\put(5,2){\Line(-1,-1)}
\put(3,3){\Line(-2,-1)}
\put(3,3){\Line(2,-1)}
\put(3,3){\circle*{1}}
\end{picture}}}= \EE_{{\begin{picture}(8,7)
\put(1,3){\circle*{0.7}}
\put(2,4){\circle*{0.7}}
\put(3,3){\circle*{0.7}}
\put(2,4){\Line(-1,-1)}
\put(2,4){\Line(1,-1)}
\put(4,5){\circle*{0.7}}
\put(5,4){\circle*{0.7}}
\put(4,5){\Line(-2,-1)}
\put(4,5){\Line(1,-1)}
\put(6,6){\circle*{0.7}}
\put(7,4){\circle*{0.7}}
\put(8,5){\circle*{0.7}}
\put(8,5){\Line(-1,-1)}
\put(6,6){\Line(-2,-1)}
\put(6,6){\Line(2,-1)}
\put(6,6){\circle*{1}}
\end{picture}}}
\ \ .
\end{equation}

\begin{equation}
\EE_{43512} \EE_{4312} = \EE_{984351267} :\qquad
\EE_{{\begin{picture}(5,5)
\put(1,3){\circle*{0.7}}
\put(2,4){\circle*{0.7}}
\put(3,2){\circle*{0.7}}
\put(4,1){\circle*{0.7}}
\put(3,2){\Line(1,-1)}
\put(5,3){\circle*{0.7}}
\put(5,3){\Line(-2,-1)}
\put(2,4){\Line(-1,-1)}
\put(2,4){\Line(3,-1)}
\put(2,4){\circle*{1}}
\end{picture}}} \EE_{{\begin{picture}(4,4)
\put(1,2){\circle*{0.7}}
\put(2,3){\circle*{0.7}}
\put(3,2){\circle*{0.7}}
\put(4,1){\circle*{0.7}}
\put(3,2){\Line(1,-1)}
\put(2,3){\Line(-1,-1)}
\put(2,3){\Line(1,-1)}
\put(2,3){\circle*{1}}
\end{picture}}}= \EE_{{\begin{picture}(9,8)
\put(1,4){\circle*{0.7}}
\put(2,5){\circle*{0.7}}
\put(3,3){\circle*{0.7}}
\put(4,2){\circle*{0.7}}
\put(3,3){\Line(1,-1)}
\put(5,4){\circle*{0.7}}
\put(5,4){\Line(-2,-1)}
\put(2,5){\Line(-1,-1)}
\put(2,5){\Line(3,-1)}
\put(6,6){\circle*{0.7}}
\put(6,6){\Line(-4,-1)}
\put(7,7){\circle*{0.7}}
\put(8,6){\circle*{0.7}}
\put(9,5){\circle*{0.7}}
\put(8,6){\Line(1,-1)}
\put(7,7){\Line(-1,-1)}
\put(7,7){\Line(1,-1)}
\put(7,7){\circle*{1}}
\end{picture}}}
\ \ .
\end{equation}

Notice that this operation coincides with the \emph{under} operation defined
by Loday-Ronco in~\cite{LR1}.

\begin{note}{\rm
The basis change from $\EE$ to $\HH$ has interesting properties:
it sends the $\EE$ trees with only left edges to their reversed tree on
$\HH$. Moreover, any $\EE_T$ is an alternating sum of $2^k$ $\HH$ functions,
where $k$ is the size of $T$ minus the length of the left edge starting from
the root of $T$.
}
\end{note}

\subsection{$\PBT$ and $\PBT^*$ as free algebras and isomorphic Hopf algebras}
\label{secisoms}

As a consequence of the existence of multiplicative bases on $\PBT$ with a
very simple product, $\PBT$ is free as an algebra (it is the algebra of a free
monoid).

Before stating and proving this result, let us recall that a permutation
$\sigma$ is \emph{connected} if it cannot be written as a shifted
concatenation $\sigma=u\sconc v$, and \emph{anticonnected} if its mirror
image $\overline{\sigma}$ is connected.

\begin{thm}
\label{PBTlibre}
The algebra $\PBT$ is free over the $\Pp_T$'s (or the $\HH_T$'s) where $T$
runs over 
 trees whose root has no right son.
In other words, $\PBT$ is free over the $\Pp_T$'s (or the $\HH_T$'s) where $T$
runs over trees whose canonical words are anticonnected permutations.
\end{thm}

\begin{pf}
The two statements of the theorem are equivalent: by definition, a tree whose
root has no right son has an anticonnected canonical word, and conversely.
Now, thanks to Section~\ref{multiplicative}, we know that the matrix that
expresses $\Pp_T$ on the basis $\HH_T$ is triangular with 1 on the main
diagonal. Moreover, the statement is obvious on the $\HH$'s thanks to their
product formula: if $\sigma$ is anticonnected, $\HH_\sigma$ cannot be
obtained by multiplication of smaller $\HH$ elements. Conversely, if $\sigma$
is not anticonnected, it can be written as $\sigma=u[k]\cdot v$ where $u$ and
$v$ are canonical words, since canonical words are the words avoiding the
pattern $132$.
\qed
\end{pf}

\begin{note}{\rm
The same theorem is true with the $\EE$'s instead of the $\HH$'s. As already
said, the product of the $\EE$'s is the shifted concatenation of $u$ and $v$
written as $u\cdot v[k]$. This proves that $\PBT$ is free over the $\Pp_T$
where $T$ runs over trees whose root has no left son.
}
\end{note}

Let us now move to $\PBT^*$. A few checks on small examples suggest that
$\PBT^*$ also is free on anticonnected canonical words. We could apply the
same  techniques to prove that it is indeed the case but we will
proceed in another way.
If $\PBT$ and $\PBT^*$ are both
free, they are isomorphic as algebras. 
We are going to prove that they are not only isomorphic as algebras but
also as Hopf algebras that will, in particular, prove that $\PBT^*$ is a free
algebra.

Let us first define a new basis in $\PBT^*$.
Let $T$ be a tree. $\Qq'_T$ is defined as
\begin{equation}
\label{qprimetoq}
\Qq'_T := \sum_{\sigma;\,\,\pp(\sigma)=T} \Qq_{\pp(\sigma^{-1})}\,.
\end{equation}

For example, representing the trees by their canonical words:
\begin{equation}
\Qq'_{231} = \Qq_{312} \quad;\quad\Qq'_{312} = \Qq_{231} + \Qq_{312}\,,
\end{equation}
\begin{equation}
\Qq'_{4213} = \Qq_{3241} + \Qq_{3412} + \Qq_{4213}\,,
\end{equation}
\begin{equation}
\Qq'_{54213} = \Qq_{43521} + \Qq_{45231} + \Qq_{45312} + \Qq_{53241} +
\Qq_{53412} + \Qq_{54213}\,,
\end{equation}
\begin{equation}
\Qq'_{53412} = \Qq_{45231} + \Qq_{45312} + \Qq_{52341} + \Qq_{53241} +
2\,\Qq_{53412} + \Qq_{54123} + \Qq_{54213}\,.
\end{equation}

As one can observe on these  examples, the smallest canonical word
in the expression of $\Qq'_\sigma$  as a sum of $\Qq_\tau$ is $\sigma^{-1}$.
This result shows that $\Qq'_\sigma$ is a basis of $\PBT^*$.

\begin{thm}
\label{qprimes}
The set $\Qq'_T$ where $T$ runs over the set of planar binary trees is a basis
of $\PBT^*$. Moreover, the matrix $M_{\Qq',\Qq}$ is triangular for the right
order:
\begin{equation}
\Qq'_\sigma=\Qq_{\sigma^{-1}} + \sum\Qq_\tau\,,
\end{equation}
where $\tau$ runs over some set of  canonical words greater than $\sigma^{-1}$
for the lexicographic order.
\end{thm}

Before completing the proof of the theorem, let us mention a simple but useful
lemma:

\begin{lem}
\label{dtsaills}
Let $\sigma$ be a permutation and $T$ its decreasing tree. Consider the
sequence of right sons starting from the root. The number of such right sons
is given by the saillances of $\sigma$. Moreover, the length of the left
subtrees attached to each son starting from the bottom-most one is given by
the saillances sequence of $\sigma$.
\end{lem}

\begin{pf}
The statement of the theorem is equivalent to the fact that, for any
permutation, the inverse of its canonical word is smaller than or equal to
the canonical word of its inverse, with equality iff the permutation is a
canonical word. Translating these facts with decreasing trees leads to the
following equivalent formulation: for any permutation $\sigma$, the canonical
word of the unlabeled shape of the decreasing tree of its canonical word is
smaller than or equal to the canonical word of the unlabeled shape of its
decreasing tree with equality iff $\sigma$ is a canonical word.
Let us prove this result.

\smallskip
First, notice that if $\sigma$ is of size $n$ and ends with an $n$, the result
is equivalent to the same statement for $\sigma'$ obtained by removing $n$
from $\sigma$. So we can assume that $\sigma$ does not end with $n$.
Moreover, it is obvious that the saillance sequence of $\sigma$ is greater
than or equal to the saillance sequence of its canonical word.
So, if the saillance sequence of $\sigma$ is different from the one of its
canonical word, thanks to Lemma~\ref{dtsaills}, its decreasing tree is
strictly greater than the decreasing tree of its canonical word.
If it is not the case, one restricts the permutation and its canonical word
to each interval between two saillances and iterate. On the trees, this
operation consists in computing the left subtrees associated with the right
sons of the root starting from the bottom-most one. By induction, this proves
the theorem.
\qed
\end{pf}

Let us now define a linear map $\phi$  from $\PBT$ to $\PBT^*$ by
\begin{equation}
\phi(\Pp_T) := \Qq'_T.
\end{equation}

\begin{thm}
\label{isomPBTavecdual}
The map $\phi$ induces a Hopf algebra isomorphism from $\PBT$ to
$\PBT^*$. In other words, one has:
\begin{equation}
\label{prodphi}
\phi(\Pp_{T'}\Pp_{T''}) = \phi(\Pp_{T'}) \phi(\Pp_{T''}) =
\Qq_{T'}' \Qq_{T''}''\,,
\end{equation}

\begin{equation}
\label{coprodphi}
(\phi\otimes\phi)(\Delta\Pp_{T}) = \Delta\phi(\Pp_{T}) = \Delta\Qq_{T}'\,.
\end{equation}
\end{thm}

\begin{pf}
First, $\phi$ is a bijection, since $\Qq'_T$ is a basis of $\PBT^*$.
Now, $\phi$ is a composition of Hopf  morphisms: it consists in the embedding
of $\PBT$ in $\FQSym$ composed with the morphism that sends $\F_\sigma$ to
$\G_{\sigma^{-1}}$ then composed with the morphism that sends $\G_\sigma$ to
its equivalence sylvester class in $\PBT^*$.
So $\phi$ is a Hopf isomorphism and both Equations~(\ref{prodphi})
and~(\ref{coprodphi}) hold.
\qed
\end{pf}

As a corollary, $\PBT$ and $\PBT^*$ are isomorphic as algebras. Since $\PBT$
is a free algebra, the same is true of $\PBT^*$.

\begin{cor}
\label{PBTetlibre}
The algebra $\PBT^*$ is free over the functions $\Qq_T$ (and $\Qq_T'$) where
$T$ runs over trees whose root has no right son.
The algebra $\PBT^*$ is free over the functions $\Qq_T$ (and $\Qq_T'$) where
$T$ runs over trees whose root has no left son.
\end{cor}

\begin{note}{\rm
One can use the isomorphism to build multiplicative bases of $\PBT^*$: they
are the sums of $\Qq'_T$ functions over upper or lower intervals of the Tamari
lattice.
}
\end{note}

\subsection{Primitive elements}

It is well-known that the dual basis of a multiplicative basis restricted to
indecomposable elements, is a basis of the Lie algebra of primitive elements
of the dual.
Since we have two mutiplicative bases on the $\PBT$ side, we then obtain two
different bases of primitive elements on $\PBT^*$. We could have worked out
the multiplicative bases on the $\PBT^*$ side but this would have been useless
since we have an explicit isomorphism of $\PBT$ to $\PBT^*$.
We obtain in this way a description of the primitive elements which differs
from that of~\cite{Ro}.

Let us denote by $\Prima_T$ (resp. $\Primb_T$) the dual bases of the $\HH_T$
(resp. $\EE_T$). The basis $\Prima$ in an analog of the basis of monomial
symmetric functions, whereas the basis $\Primb$ is an analog of the
forgotten symmetric functions. The following results hold:

\begin{thm}
The Lie algebra of primitive elements of $\PBT^*$ is spanned by the
$\Prima_T$'s where $T$ runs over trees whose roots have no right son.
The Lie algebra of primitive elements of $\PBT^*$ is spanned by the
$\Primb_T$'s where $T$ runs over trees whose roots have no left son.
\end{thm}

The first matrices $M_{\Prima_T,\Qq_T}$ and $M_{\Primb_T,\Qq_T}$ for trees up
to $4$ is respectively given in Figures~\ref{m2q23}--\ref{fm2q4}.

\begin{note}{\rm
Since $M_{\Prima_T,\Qq_T}$ is the transpose of $M_{\Pp_T,\HH_T}$, the
expression of $\Prima_T$ on $\Qq_T$ is derived from the M\"obius inversion of
the Tamari lattice.
}
\end{note}

\subsection{Embeddings and quotients}
\label{morphismes}

\subsubsection{The full diagram of embeddings}

In~\cite{LR1}, Loday and Ronco defined different morphisms starting from or
getting to $\PBT$. These morphisms can be naturally understood and realized in
our framework since we have non-commutative polynomial realizations of all of
those: $\FQSym$, $\PBT$, and $\NCSF$. Indeed, the morphisms become trivial:
all algebras are included in the same non-commutative polynomial algebra.

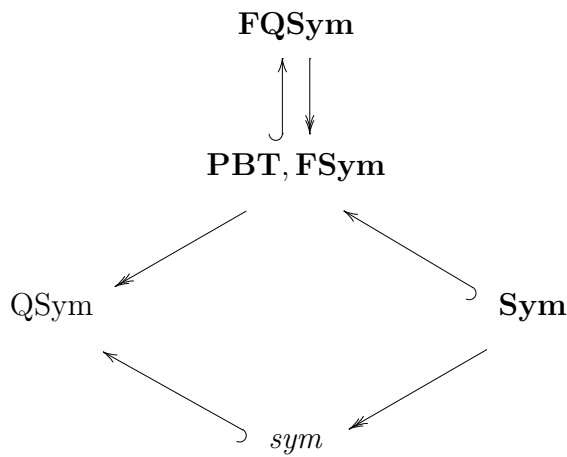
\begin{figure}[ht]
\begin{equation*}
\xymatrix@R=1cm@C=1cm{
          &  *+<5mm>{\FQSym}\ar@<1ex>@{>>}[d]       & \\
          &  *+<5mm>{\PBT, \FSym}\ar@<1ex>@{>>}[dl]\ar@<1ex>@{^{(}->}[u]       & \\
*+<5mm>{\QSym} &             & *+<5mm>{\NCSF}\ar@<1ex>@{>>}[dl]\ar@<1ex>@{^{(}->}[ul]       & \\
          &  *+<5mm>{\sym}\ar@<1ex>@{^{(}->}[ul]  & \\
}
\end{equation*}
\caption{\label{diagalg}Morphisms between known related Hopf algebras.}
\end{figure}

Let us first present the general diagram containing all these algebras and a
few other ones (see Figure~\ref{diagalg}).

The algebra $\sym$ is the usual algebra of commutative Symmetric Functions.
As one can see on Figure~\ref{diagalg}, all these algebras are subalgebras or
quotients of $\FQSym$. As this algebra can be realized in the free algebra on
an infinite alphabet, it is the same for all the other algebras.
Then the up arrows just are inclusions and the down arrows are induced by
commutation rules amoung the letters of the alphabet.

Let us describe in more detail two arrows:
\begin{itemize}
\item $\NCSF$ to $\PBT$:
the algebra of non-commutative symmetric functions is generated by the
homogeneous symmetric functions $S_n$ which can be realized as the sum
of all non-decreasing words of length $n$. As a polynomial, $S_n$ is equal to
the $\Pp_T=\HH_T$ function where $T$ is the tree with $n-1$ left edges.
Since $S_n$ and $\HH_T$ are both multiplicative bases, it is obvious that one
gets a realization of $\NCSF$ inside $\PBT$ as the subalgebra generated by the
$\HH_T$ where $T$ runs over the set of right comb trees.
Moreover, the basis of elementary non-commutative symmetric functions
$\Lambda_n$ is realized as the sum of all decreasing words of size $n$ that
happens to be $\EE_T=\Pp_T$, where $T$ is the tree with $n-1$ right edges.
Thus, the linear basis of $\Lambda^I$ is realized as the $\EE_T$ where $T$
runs over the left comb trees.
\item $\PBT^*$ to $\QSym$: 
the algebra $\PBT^*$ is the specialization of $\FQSym^*$ to a ``sylvester
alphabet'', an ordered alphabet satisfying the sylvester relations whereas
$\QSym$ is the specialization of both $\PBT^*$ and $\FQSym^*$ to a commutative
alphabet.
\end{itemize}

\begin{note}{\rm
Figure~\ref{diagalg} is a commutative diagram. Indeed, both compound morphisms
$\NCSF\to\sym\to\QSym$ and $\NCSF\to\PBT\to\FQSym\to\FQSym^*\to\PBT^*\to\QSym$
amount to compute commutative images of polynomials.
}
\end{note}

\begin{note}{\rm
The pair of Hopf algebras in duality $\FSym$ and $\FSym^*$ play the same
role as $\PBT$ and $\PBT^*$.
However, there exists a difference between these algebras: the compound
morphism $\PBT\to\FQSym\to\PBT^*$ is a Hopf isomorphim whereas the compound
map $\FSym\to\FQSym\to\FSym^*$ is not an isomorphism, since it is not even
injective.
}
\end{note}

%
\section{Representation theory}
%

As mentioned in the introduction, it is known that the integers
\begin{equation}
c_{I,J} = \< R_I, R_J \> \quad |I|=|J|=n\,,
\end{equation}
where $R_I$ stands for the ribbon-Schur function of shape $I$, can be
interpreted as Cartan invariants of the $0$-Hecke algebra: the coefficient
$c_{I,J}$ is equal to the multiplicity of the simple module $S_I$ in the
indecomposable projective module $P_J$ (see~\cite{NCSF4}).

The analogy between ribbon Schur functions and the natural basis $\Pp_T$
of $\PBT$ allows one to wonder whether one can
interpret in the same way the integers 
\begin{equation}
c_{T,U} = \< \Pp_T, \Pp_U \> 
= \operatorname{Card}\{\sigma;\ \pp(\sigma)=T,\ \pp(\sigma^{-1})=U\}.
\end{equation}

Let $M^{(n)}$ be the matrix of the $c_{T,U}$ ordered in rows and columns by
the lexicographic order of the canonical words. It suffices to compute
$M^{(3)}$ (see Figure~\ref{matm3}) to understand that this matrix cannot
be a matrix of Cartan invariants: it has a $0$ on the diagonal.

\begin{figure}[ht]
\begin{equation}
\begin{pmatrix}
\,1\, &       &       &       &       \\
      & \,1\, &       &       &       \\
      &       &       & \,1\, &       \\
      &       & \,1\, & \,1\, &       \\
      &       &       &       & \,1\, \\
\end{pmatrix}
\qquad\qquad
\begin{pmatrix}
123 &     &     &     &     \\
    & 213 &     &     &     \\
    &     &     & 231 &     \\
    &     & 312 & 132 &     \\
    &     &     &     & 321 
\end{pmatrix}
\end{equation}
\caption{\label{matm3}The matrix $M^{(3)}$ and the corresponding
permutations.}
\end{figure}

Indeed, if one assumes that $\<\Pp_T,\Pp_U\>= \dim\hom_{A_n}(\Pp_T,\Pp_U)$,
or, equivalently if one assumes that the simple modules are indexed in such a
way that $S_T=P_T/{\hbox{Rad $P_T$}}$, each diagonal entry is at least $1$
since $\hom_{A_n}(\Pp_T,\Pp_T)$ contains at least the identity map.

\subsection{Combinatorial analysis of the scalar product}

\subsubsection{The Gram matrices}

However, the Gram matrices $M^{(n)}$ have an interesting block structure. This
leads to enquire whether there exists a simple transformation building a more
interesting sequence of matrices.
We already solved this question in Section~\ref{secisoms} and more precisely
in Theorem~\ref{qprimes}. Indeed, if one orders the rows of $M^{(n)}$ with the
lexicographic order of their canonical words and the columns with the
lexicographic order of the inverses of the canonical words, the matrix $M$
becomes the matrix expressing the $\Qq'_T$ on the $\Qq_T$ basis.

We will now present the block structure in order to get the right order
on trees and its interpretation in terms of the scalar product inherited from
$\FQSym$.

\smallskip
Let $T$ be a planar binary tree. The \emph{skeleton} of $T$ is the pair of
integers $(k,l)$ defined by
\begin{itemize}
\item $k\leq n-1$ is the greatest integer such that
$w_T(n)=n,\ldots,w_T(n-k+1)=n-k+1$, say, the number of fixed points at the end
of $w_T$ (minus $1$ if $w_T$ is the identity permutation),
\item $l$ is the number of saillances of $\sigma_T$, after one has removed
its last $k$ elements.
\end{itemize}

One can geometrically define the skeleton of a tree $T$ as the part of
$T$ composed of the highest sequence of right sons and of vertices greater
than this sequence. Figure~\ref{squelette} shows the skeleton of
$(8,9,7,5,4,6,1,2,3,10,11)$. The number of fixed points at the end of the
permutation (two in the example) corresponds to the size of the left edge
minus $1$, and the number of saillances of the remaining permutation (four in
the example) corresponds to the size of the right edge.

\begin{figure}[ht]
\centerline{
\begingroup
\setlength{\unitlength}{2000sp}%
%
\begin{picture}(7503,3826)(2287,-4074)
{\thinlines
\put(4801,-361){\circle{212}}
}%
{\put(4201,-961){\circle{212}}
}%
{\put(3601,-1561){\circle{212}}
}%
{\put(3001,-2161){\circle{212}}
}%
{\put(2401,-2761){\circle{212}}
}%
{\put(4201,-2161){\circle{212}}
}%
{\put(4801,-2761){\circle{212}}
}%
{\put(5401,-3361){\circle{212}}
}%
{\put(4801,-3961){\circle{212}}
}%
{\put(4201,-3361){\circle{212}}
}%
{\put(3601,-2761){\circle{212}}
}%
\thicklines
{\put(4726,-436){\line(-1,-1){450}}
}%
{\put(4126,-1036){\line(-1,-1){450}}
}%
{\put(3526,-1636){\line(-1,-1){450}}
}%
{\put(2926,-2236){\line(-1,-1){450}}
}%
{\put(4126,-2236){\line(-1,-1){450}}
}%
{\put(5326,-3436){\line(-1,-1){450}}
}%
{\put(3676,-1636){\line( 1,-1){450}}
}%
{\put(4276,-2236){\line( 1,-1){450}}
}%
{\put(4876,-2836){\line( 1,-1){450}}
}%
{\put(3675,-2835){\line( 1,-1){450}}
}%
{\thinlines
\put(9076,-361){\circle{212}}
}%
{\put(8476,-961){\circle{212}}
}%
{\put(7876,-1561){\circle{212}}
}%
{\put(8476,-2161){\circle{212}}
}%
{\put(9076,-2761){\circle{212}}
}%
{\put(9676,-3361){\circle{212}}
}%
\thicklines
{\put(9001,-436){\line(-1,-1){450}}
}%
{\put(7951,-1636){\line( 1,-1){450}}
}%
{\put(8551,-2236){\line( 1,-1){450}}
}%
{\put(9151,-2836){\line( 1,-1){450}}
}%
\end{picture}
\endgroup
}
\caption{\label{squelette}A tree and its skeleton.}
\end{figure}

We can now describe the block structure of $M^{(n)}$: let us say that two
trees $T$ and $U$ are \emph{in the same block} if there exists a power of
$M^{(n)}$ in which the coefficient $(T,U)$ is nonzero. Notice that this
relation is symmetric since $M^{(n)}$ is itself symmetric.

\begin{thm}
\label{blockthm}
Two trees are in the same block iff they have the same size and the same
skeleton.
\end{thm}

To prove this property, we will need a definition and two simple lemmas:

\begin{defn}
Let $\sigma$ and $\sigma'$ be two permutations such that there exists three
indices $i<j<k$ such that the restriction of both permutations to these
indices are two words $acb$ and $bca$ with $a<b<c$.
We say that $\sigma$ and $\sigma'$ are \emph{co-sylvester adjancent}.
This allows us to define the co-sylvester equivalence by transitive closure of
the co-sylvester adjacence.
\end{defn}

\begin{lem}
Two permutations $\sigma$ and $\sigma'$ are co-sylvester-adjacent
(respectively co-sylvester-equivalent) iff $\sigma^{-1}$ and $\sigma'^{-1}$
are sylvester-adjacent (resp. sylvester-congruent).
\end{lem}

\begin{lem}
The greatest word for the lexicographic order associated with a tree of size
$n$ of skeleton $(k,l)$ is given by
\begin{equation}
w_{k,l,n}:=(n-k)\cdots(n-k-l+2)(n-k-l)\cdots1.(n-k-l+1).(n-k+1)\cdots n.
\end{equation}
\end{lem}

For example, $w_{2,4,9}=765321489$. There exists another description of $w$:
given the skeleton, build a tree by attaching the tree of the correct size
only composed of right edges to the left of the left-most node of the
skeleton. Then $w$ is the canonical word associated with this tree. It is a
consequence of Lemma~\ref{dtsaills}.

Let us now prove Theorem~\ref{blockthm}.

\begin{pf}
By definition, two trees are in the same matrix block iff there is a path
going from the first one to the second one consisting of pairs of trees, the
first one being the binary tree of a permutation and the second one
being the binary tree of its inverse.
It is obvious that both binary trees have same skeleton since
$\sigma$ and $\sigma^{-1}$ have the same number of fixed points at the end and
the same number of saillances after having removed the previous fixed points.

Conversely, let us prove that all trees with a given skeleton are connected.
Thanks to Theorem~\ref{qprimes}, the permutations corresponding to a given row
of $M^{(n)}$ are sylvester-congruent whereas the permutations corresponding to
a given column are co-sylvester-equivalent. Consider a skeleton $(k,l)$. Then
proving that all trees are connected is equivalent to prove that all
permutations of size $n$ having skeleton $(k,l)$ are connected using both
sylvester and co-sylvester relations. This will come from the fact that
they all are connected to $w_{k,l,n}$.

By induction, we can restrict to the case $k=0$ and to canonical sylvester
permutations. Let $\sigma$ be a canonical permutation of size $n$ and skeleton
$(0,l)$. If $\sigma$ does not begin with $n$ then exchange the two neighboors
of $n$ by a co-sylvester rewriting and take the sylvester rewriting on these
three elements. This permutation is greater than $\sigma$ so is its canonical
word. Iterating this process, one ends with a permutation beginning with $n$.
If $\sigma$ begins with $n$ then by induction on $n$, it is connected to
$w_{k,l,n}$. So all permutations of size $n$ and skeleton $(k,l)$ are
connected to $w_{k,l,n}$.
\qed
\end{pf}

The proof of the next theorem directly follows from Theorem~\ref{qprimes}.

\begin{thm}
Let $\nu$ be the involution on trees defined as $\nu(T)=T'$, where
$T'=\pp(w_T^{-1})$. Then
\begin{itemize}
\item the involution $\nu$ preserves the blocks: ($T$ and $\nu(T)$ have same
skeleton),
\item the matrix $C^{(n)}$ defined by
\begin{equation}
C^{(n)}(T,U) = \< \Pp_{T}, \Pp_{\nu(U)} \>
\end{equation}
is block lower unitriangular if one orders the trees, first by skeleton, then
by lexicographic order on the canonical words of each skeleton class of
trees.
\end{itemize}
\end{thm}

Figure~\ref{matc} contains the first matrices $C^{(n)}$, skipping the zero
entries to allow instantaneous reading.
The order of the trees in rows and columns corresponds to the lexicographic
order on their canonical words:
\begin{itemize}
\item $12$; $21$,
\item $123$; $213$; $231$, $312$; $321$,
\item $1234$; $2134$; $2314$, $3124$; $3214$; $2341$, $3241$, $3412$, $4123$,
$4213$; $3421$, $4231$, $4312$; $4321$.
\end{itemize}
\begin{figure}[ht]
\scriptsize
\begin{center}
$\begin{pmatrix}
1& \\
 &1
\end{pmatrix}$
\quad
$\begin{pmatrix}
1& & & & \\
 &1& & & \\
 & &1& & \\
 & &1&1& \\
 & & & &1
\end{pmatrix}$
\quad
$\left(\begin{array}{cccccccccccccc}
1& & & & & & & & & & & & & \\
 &1& & & & & & & & & & & & \\
 & &1& & & & & & & & & & & \\
 & &1&1& & & & & & & & & & \\
 & & & &1& & & & & & & & & \\
 & & & & &1& & & & & & & & \\
 & & & & & &1& & & & & & & \\
 & & & & &1&1&1& & & & & & \\
 & & & & &1& &1&1& & & & & \\
 & & & & & &1&1& &1& & & & \\
 & & & & & & & & & &1& & & \\
 & & & & & & & & & &1&1& & \\
 & & & & & & & & & &1&1&1& \\
 & & & & & & & & & & & & &1
\end{array}\right)$
\end{center}
\caption{\label{matc}The matrices $C^{(n)}$ for $n=2$, $3$, $4$.}
\end{figure}

\subsubsection{Combinatorics of the Gram matrices}

Let us now study more precisely the block structure of our matrices.
We first need a few classical definitions.
Define the \emph{Catalan triangle} (see~\cite{Slo}) (resp. the \emph{first
kind Stirling triangle}) as the triangular matrix $A$ (resp. $B$) whose
coefficient $a_{i,j}$ (resp. $b_{i,j}$) is the coefficient of $t^i u^j$ in the
respective expressions:

\begin{equation}
 \sum_{n=1}^\infty{t^n\ \sum_{k=0}^{n-1}{\frac{n-k}{n+k}\binom{n+k}{k}u^k}},
\quad\quad
\sum_{n=1}^\infty{t^n \prod_{k=1}^{n-1}{(1+k\ u)}}.
\end{equation}

The well-known combinatorial interpretation of these numbers is the following:
$a_{i,j}$ is the number of planar binary trees of size $i$ whose number of
elements not belonging to the sequence of right sons starting from the root is
$j-1$.
The coefficient $b_{i,j}$ is the number of permutations of size $i$ whose
saillances number is $i-j+1$.
Both triangles are represented Figure~\ref{catsti} where we put into
parenthesis the trees having no right son at the root, and the permutations
having only one saillance.
Notice that the same Catalan triangle has been encountered by
Aval-Bergeron-Bergeron in~\cite{ABB} while studying the quotient of the
algebra of polynomials by the ideal generated by quasi-symmetric polynomials
without constant term. The relation between both constructions remains
mysterious.

\begin{figure}[ht]
\centerline{$\begin{array}{cccccccc}
1 &     &     &     &      &       &       \\
1 & (1) &     &     &      &       &       \\
1 &   2 & (2) &     &      &       &       \\
1 &   3 &   5 & (5) &      &       &       \\
1 &   4 &   9 &  14 & (14) &       &       \\
1 &   5 &  14 &  28 &   42 & (42)  &       \\
1 &   6 &  20 &  48 &   90 & 132   & (132)
\end{array}$
\qquad\qquad
$\begin{array}{cccccccc}
1 &     &     &     &      &       &       \\
1 & (1) &     &     &      &       &       \\
1 &   3 & (2) &     &      &       &       \\
1 &   6 &  11 & (6) &      &       &       \\
1 &  10 &  35 &  50 & (24) &       &       \\
1 &  15 &  85 & 225 &  274 & (120) &       \\
1 &  21 & 175 & 735 & 1624 & 1764  & (720) \\
\end{array}$
}
\caption{\label{catsti}The Catalan and first kind Stirling triangles.}
\end{figure}

\begin{thm}
Let us consider a block of skeleton $(k,l)$ of $C^{(n)}$ for a given $n$.

\begin{itemize}
\item The number $m_n(k,l)$ of rows of this block is given by the
$(n-k-l+1)$-th number of the $(n-k)$-th row of the Catalan triangle.
\item The sum $d_n(k,l)$ of the entries of this block is given by the
$(n-k-l-1)$-th number of the $(n-k)$-th row of the first kind Stirling
triangle.
\end{itemize}
\end{thm}

\subsection{A conjectural representation theoretical interpretation}

\begin{defn}
A \emph{tower of algebras} is a pair
$\left((A_n)_{n\in\NN},\, (\rho_{i,j})_{i,j\in\NN}\right)$
where the $A_n$'s are algebras, and for all $i,j$ the map $\rho_{i,j}$ is
an algebra embedding of $A_i\tensor A_j$ into $A_{i+j}$ such that
\begin{equation}\label{eq.assocRho}
\rho_{i+j,k} \circ  (\rho_{i,j} \tensor \Id_{A_k}) =
\rho_{i,j+k} \circ  (\Id_{A_i} \tensor \rho_{j,k})\,.
\end{equation}
\end{defn}

Notice that Equation~(\ref{eq.assocRho}) amounts to require that the direct
sum of the maps $(x,y) \mapsto \rho_{i,j}(x\otimes y)$ if $x\in A_i$ and $y\in
A_j$ defines an associative product on the direct sum $\bigoplus A_i$ which is
compatible with the structure of the $A_i$.

For any tower of algebras, the induction process with respect to the
embeddings $\rho_{i,j}$ defines an algebra structure on the direct sums of
the Grothendieck groups
\begin{equation}
{\mathcal G}=\bigoplus_{n\ge 0}G_0(A_n)\,,\quad
{\mathcal K}=\bigoplus_{n\ge 0}K_0(A_n).
\end{equation}
Similarly, the restriction defines a coproduct in such a way that
${\mathcal G}$ and ${\mathcal K}$ are equipped with two mutually dual Hopf
algebra structures.

An example of a tower of algebras is the tower of algebras of the symmetric
groups together with the linear maps extending the group inclusions
$\SG_i\times\SG_j \mapsto \SG_{i+j}$. It is well-known that this leads to the
self-dual Hopf algebra of symmetric functions (see~\cite{G}). Replacing the
symmetric group by its degenerated Hecke algebra leads to the dual pair
$(\QSym,\Sym)$. And it is likely that the pair $(\PBT,\PBT^*)$ comes from a
similar construction.

For any tower of algebras $(A_n)$ having such Gram matrices as Cartan
invariants, the skeletons $(k,l)$ correspond to the blocks (the indecomposable
subalgebras) $B_n(k,l)$ of $A_n$. Then, $d_n(k,l)$ is the dimension of
$B_n(k,l)$, and $m_n(k,l)$ is the number of its simple modules.

Notice that there might \emph{a priori} exist many non isomorphic towers of
algebras such that the $C^{(n)}$ matrices are their Cartan invariants.

We computed the quivers corresponding to each block of the matrix for
$n\leq6$ with the constraint of providing the smallest possible amount of
arrows (equivalently the smallest possible amount of relations).
The structure of these quivers and their relations seem to have a certain
regularity but they unfortunately remain unsufficiently understood to allow
us to describe these algebras for any~$n$.
Nevertheless, we conjecture the following result:

\begin{conj}
There exists a tower of algebras $(A_n)$ such that the $C^{(n)}$ are their
matrices of Cartan invariants.
\end{conj}
 
In particular, one should have $\dim A_n=n!$.
We can also propose a more precise conjecture.


\begin{conj}
There exists a tower of algebras $A_n$, with a basis
$(e_\sigma)_{\sigma\in\SG_n}$ such that:

\begin{itemize}
\item the restriction to canonical words of the morphism $\rho_{m,n}$ is given
by the product of the corresponding $\Pp_T$ functions.
In this setting, the indecomposable projective modules of $A_n$ are left
ideals $P_T=A_n e_{\sigma_T}$, and therefore are in bijection with the planar
binary trees of size $n$.
\item If one endows ${\mathcal K}$ with the induction product
$[M]\cdot [N]=[M\otimes_\CC N\uparrow_{A_m\otimes A_n}^{A_{m+n}}]$,
the map ${\mathcal K}\rightarrow\PBT$ sending the class of the module
$P_T$ on the polynomial $\Pp_T$ is a ring isomorphism.
\end{itemize}
\end{conj}

%
\section{Conclusion}
%

Since its discovery in the mid-seventies, the plactic monoid has, for a long
time, been considered as a very singular object. It needed the discovery of
quantum groups (independently due to Drinfeld and Jimbo about 1985)
(see~\cite{Jim}), and Kashiwara's theory of crystal bases (1991)
(see~\cite{Kashiwara91}), to discover the plactic monoids associated with all
semi-simple Lie algebras (see~\cite{LascouxAl95,Littelmann96}). But even this
point of view does not tell everything about plactic monoids. The hypoplactic
monoid (see~\cite{NCSF4}), which is to quasi-symmetric functions what the
ordinary plactic monoid is to ordinary symmetric functions, was obtained from
a non standard version of the quantum linear group, and is not taken into
account by the theory of crystal bases. This raises a first question, to find
a quantum group interpretation of the sylvester monoid, and a second one, to
characterize and classify all similar monoids.

%
%
\newpage
\section{Tables}
%
%

In this Section, we give the transition matrices between various bases in
degree $n\leq4$.
Rows and columns of those matrices correspond to binary trees on $n$ nodes
arranged as follows:

\begin{figure}[ht]
\begin{equation}
\left[{\begin{picture}(2,3)
\put(1,1){\circle*{0.7}}
\put(2,2){\circle*{0.7}}
\put(2,2){\Line(-1,-1)}
\put(2,2){\circle*{1}}
\end{picture}},{\begin{picture}(2,3)
\put(1,2){\circle*{0.7}}
\put(2,1){\circle*{0.7}}
\put(1,2){\Line(1,-1)}
\put(1,2){\circle*{1}}
\end{picture}}\ \right]
\end{equation}
\caption{\label{ord2}Order on trees of size $2$.}
\end{figure}

\begin{figure}[ht]
\begin{equation}
\left[{\begin{picture}(3,4)
\put(1,1){\circle*{0.7}}
\put(2,2){\circle*{0.7}}
\put(2,2){\Line(-1,-1)}
\put(3,3){\circle*{0.7}}
\put(3,3){\Line(-1,-1)}
\put(3,3){\circle*{1}}
\end{picture}},{\begin{picture}(3,4)
\put(1,2){\circle*{0.7}}
\put(2,1){\circle*{0.7}}
\put(1,2){\Line(1,-1)}
\put(3,3){\circle*{0.7}}
\put(3,3){\Line(-2,-1)}
\put(3,3){\circle*{1}}
\end{picture}},{\begin{picture}(3,4)
\put(1,3){\circle*{0.7}}
\put(2,1){\circle*{0.7}}
\put(3,2){\circle*{0.7}}
\put(3,2){\Line(-1,-1)}
\put(1,3){\Line(2,-1)}
\put(1,3){\circle*{1}}
\end{picture}},{\begin{picture}(3,3)
\put(1,1){\circle*{0.7}}
\put(2,2){\circle*{0.7}}
\put(3,1){\circle*{0.7}}
\put(2,2){\Line(-1,-1)}
\put(2,2){\Line(1,-1)}
\put(2,2){\circle*{1}}
\end{picture}},{\begin{picture}(3,4)
\put(1,3){\circle*{0.7}}
\put(2,2){\circle*{0.7}}
\put(3,1){\circle*{0.7}}
\put(2,2){\Line(1,-1)}
\put(1,3){\Line(1,-1)}
\put(1,3){\circle*{1}}
\end{picture}}\ \right]
\end{equation}
\caption{\label{ord3}Order on trees of size $3$.}
\end{figure}

\begin{figure}[ht]
\begin{equation}
\left[
\begin{split}
&{\begin{picture}(4,5)
\put(1,1){\circle*{0.7}}
\put(2,2){\circle*{0.7}}
\put(2,2){\Line(-1,-1)}
\put(3,3){\circle*{0.7}}
\put(3,3){\Line(-1,-1)}
\put(4,4){\circle*{0.7}}
\put(4,4){\Line(-1,-1)}
\put(4,4){\circle*{1}}
\end{picture}},{\begin{picture}(4,5)
\put(1,2){\circle*{0.7}}
\put(2,1){\circle*{0.7}}
\put(1,2){\Line(1,-1)}
\put(3,3){\circle*{0.7}}
\put(3,3){\Line(-2,-1)}
\put(4,4){\circle*{0.7}}
\put(4,4){\Line(-1,-1)}
\put(4,4){\circle*{1}}
\end{picture}},{\begin{picture}(4,5)
\put(1,3){\circle*{0.7}}
\put(2,1){\circle*{0.7}}
\put(3,2){\circle*{0.7}}
\put(3,2){\Line(-1,-1)}
\put(1,3){\Line(2,-1)}
\put(4,4){\circle*{0.7}}
\put(4,4){\Line(-3,-1)}
\put(4,4){\circle*{1}}
\end{picture}},{\begin{picture}(4,5)
\put(1,4){\circle*{0.7}}
\put(2,1){\circle*{0.7}}
\put(3,2){\circle*{0.7}}
\put(3,2){\Line(-1,-1)}
\put(4,3){\circle*{0.7}}
\put(4,3){\Line(-1,-1)}
\put(1,4){\Line(3,-1)}
\put(1,4){\circle*{1}}
\end{picture}},{\begin{picture}(4,5)
\put(1,2){\circle*{0.7}}
\put(2,3){\circle*{0.7}}
\put(3,2){\circle*{0.7}}
\put(2,3){\Line(-1,-1)}
\put(2,3){\Line(1,-1)}
\put(4,4){\circle*{0.7}}
\put(4,4){\Line(-2,-1)}
\put(4,4){\circle*{1}}
\end{picture}},{\begin{picture}(4,5)
\put(1,3){\circle*{0.7}}
\put(2,2){\circle*{0.7}}
\put(3,1){\circle*{0.7}}
\put(2,2){\Line(1,-1)}
\put(1,3){\Line(1,-1)}
\put(4,4){\circle*{0.7}}
\put(4,4){\Line(-3,-1)}
\put(4,4){\circle*{1}}
\end{picture}},{\begin{picture}(4,5)
\put(1,4){\circle*{0.7}}
\put(2,2){\circle*{0.7}}
\put(3,1){\circle*{0.7}}
\put(2,2){\Line(1,-1)}
\put(4,3){\circle*{0.7}}
\put(4,3){\Line(-2,-1)}
\put(1,4){\Line(3,-1)}
\put(1,4){\circle*{1}}
\end{picture}},\\ &{\begin{picture}(4,4)
\put(1,2){\circle*{0.7}}
\put(2,3){\circle*{0.7}}
\put(3,1){\circle*{0.7}}
\put(4,2){\circle*{0.7}}
\put(4,2){\Line(-1,-1)}
\put(2,3){\Line(-1,-1)}
\put(2,3){\Line(2,-1)}
\put(2,3){\circle*{1}}
\end{picture}},{\begin{picture}(4,5)
\put(1,4){\circle*{0.7}}
\put(2,3){\circle*{0.7}}
\put(3,1){\circle*{0.7}}
\put(4,2){\circle*{0.7}}
\put(4,2){\Line(-1,-1)}
\put(2,3){\Line(2,-1)}
\put(1,4){\Line(1,-1)}
\put(1,4){\circle*{1}}
\end{picture}},{\begin{picture}(4,4)
\put(1,1){\circle*{0.7}}
\put(2,2){\circle*{0.7}}
\put(2,2){\Line(-1,-1)}
\put(3,3){\circle*{0.7}}
\put(4,2){\circle*{0.7}}
\put(3,3){\Line(-1,-1)}
\put(3,3){\Line(1,-1)}
\put(3,3){\circle*{1}}
\end{picture}},{\begin{picture}(4,4)
\put(1,2){\circle*{0.7}}
\put(2,1){\circle*{0.7}}
\put(1,2){\Line(1,-1)}
\put(3,3){\circle*{0.7}}
\put(4,2){\circle*{0.7}}
\put(3,3){\Line(-2,-1)}
\put(3,3){\Line(1,-1)}
\put(3,3){\circle*{1}}
\end{picture}},{\begin{picture}(4,5)
\put(1,4){\circle*{0.7}}
\put(2,2){\circle*{0.7}}
\put(3,3){\circle*{0.7}}
\put(4,2){\circle*{0.7}}
\put(3,3){\Line(-1,-1)}
\put(3,3){\Line(1,-1)}
\put(1,4){\Line(2,-1)}
\put(1,4){\circle*{1}}
\end{picture}},{\begin{picture}(4,4)
\put(1,2){\circle*{0.7}}
\put(2,3){\circle*{0.7}}
\put(3,2){\circle*{0.7}}
\put(4,1){\circle*{0.7}}
\put(3,2){\Line(1,-1)}
\put(2,3){\Line(-1,-1)}
\put(2,3){\Line(1,-1)}
\put(2,3){\circle*{1}}
\end{picture}},{\begin{picture}(4,5)
\put(1,4){\circle*{0.7}}
\put(2,3){\circle*{0.7}}
\put(3,2){\circle*{0.7}}
\put(4,1){\circle*{0.7}}
\put(3,2){\Line(1,-1)}
\put(2,3){\Line(1,-1)}
\put(1,4){\Line(1,-1)}
\put(1,4){\circle*{1}}
\end{picture}}
\end{split}
\right]
\end{equation}
\caption{\label{ord4}Order on trees of size $4$.}
\end{figure}

These orders correspond to the lexicographic order on canonical words:

\begin{equation}
[12, 21] ;\quad\quad [123, 213, 231, 312, 321] ;
\end{equation}
\begin{equation}
\left[
\begin{split}
& 1234, 2134, 2314, 2341, 3124, 3214, 3241,\\
& 3412, 3421, 4123, 4213, 4231, 4312, 4321\,.
\end{split}
\right]
\end{equation}

Now, let us give the matrices $M_{\HH,\Pp}$, $M_{\EE,\Pp}$,
$M_{\EE,\HH}$, $M_{\Qq',\Qq}$, $M_{\Prima,\Qq}$, and finally $M_{\Primb,\Qq}$
for $n=2$, $3$ and $4$.

Notice that the matrix $M_{\Prima,\Qq}$ is the transpose of the inverse of
$M_{\HH,\Pp}$.
It is the same with $M_{\Primb,\Qq}$, that is the transpose of the inverse of
$M_{\EE,\Pp}$.

\newpage
~\vfill

\begin{figure}[ht]
\newcolumntype{C}{@{\hskip15pt}>{$}c<{$}}
\begin{equation*}
\left(\begin{array}{cC}
1 & 1\\
. & 1\\
\end{array}\right)
\qquad\qquad
\left(\begin{array}{cCCCC}
1 & 1 & 1 & 1 & 1\\
. & 1 & 1 & . & 1\\
. & . & 1 & . & 1\\
. & . & . & 1 & 1\\
. & . & . & . & 1\\
\end{array}\right)
\end{equation*}
\caption{\label{h2p23}The matrices $M_{\HH_T,\Pp_T}$ for $n=2,3$.}
\end{figure}

\vfill

\begin{figure}[ht]
\newcolumntype{C}{@{\hskip15pt}>{$}c<{$}}
\begin{equation*}
\left(\begin{array}{cCCCCCCCCCCCCC}
1 & 1 & 1 & 1 & 1 & 1 & 1 & 1 & 1 & 1 & 1 & 1 & 1 & 1\\
. & 1 & 1 & 1 & . & 1 & 1 & . & 1 & . & 1 & 1 & . & 1\\
. & . & 1 & 1 & . & 1 & 1 & . & 1 & . & . & 1 & . & 1\\
. & . & . & 1 & . & . & 1 & . & 1 & . & . & 1 & . & 1\\
. & . & . & . & 1 & 1 & 1 & 1 & 1 & . & . & . & 1 & 1\\
. & . & . & . & . & 1 & 1 & . & 1 & . & . & . & . & 1\\
. & . & . & . & . & . & 1 & . & 1 & . & . & . & . & 1\\
. & . & . & . & . & . & . & 1 & 1 & . & . & . & 1 & 1\\
. & . & . & . & . & . & . & . & 1 & . & . & . & . & 1\\
. & . & . & . & . & . & . & . & . & 1 & 1 & 1 & 1 & 1\\
. & . & . & . & . & . & . & . & . & . & 1 & 1 & . & 1\\
. & . & . & . & . & . & . & . & . & . & . & 1 & . & 1\\
. & . & . & . & . & . & . & . & . & . & . & . & 1 & 1\\
. & . & . & . & . & . & . & . & . & . & . & . & . & 1\\
\end{array}\right)
\end{equation*}
\caption{\label{h2p4}The matrix $M_{\HH_T,\Pp_T}$ for $n=4$.}
\end{figure}

\vfill
\newpage
~\vfill

\begin{figure}[ht]
\newcolumntype{C}{@{\hskip15pt}>{$}c<{$}}
\begin{equation*}
\left(\begin{array}{cC}
1 & .\\
1 & 1\\
\end{array}\right)
\qquad\qquad
\left(\begin{array}{cCCCC}
1 & . & . & . & .\\
1 & 1 & . & . & .\\
1 & 1 & 1 & . & .\\
1 & . & . & 1 & .\\
1 & 1 & 1 & 1 & 1\\
\end{array}\right)
\end{equation*}
\caption{\label{e2p23}The matrices $M_{\EE_T,\Pp_T}$ for $n=2,3$.}
\end{figure}

\vfill

\begin{figure}[ht]
\newcolumntype{C}{@{\hskip15pt}>{$}c<{$}}
\begin{equation*}
\left(\begin{array}{cCCCCCCCCCCCCC}
1 & . & . & . & . & . & . & . & . & . & . & . & . & .\\
1 & 1 & . & . & . & . & . & . & . & . & . & . & . & .\\
1 & 1 & 1 & . & . & . & . & . & . & . & . & . & . & .\\
1 & 1 & 1 & 1 & . & . & . & . & . & . & . & . & . & .\\
1 & . & . & . & 1 & . & . & . & . & . & . & . & . & .\\
1 & 1 & 1 & . & 1 & 1 & . & . & . & . & . & . & . & .\\
1 & 1 & 1 & 1 & 1 & 1 & 1 & . & . & . & . & . & . & .\\
1 & . & . & . & 1 & . & . & 1 & . & . & . & . & . & .\\
1 & 1 & 1 & 1 & 1 & 1 & 1 & 1 & 1 & . & . & . & . & .\\
1 & . & . & . & . & . & . & . & . & 1 & . & . & . & .\\
1 & 1 & . & . & . & . & . & . & . & 1 & 1 & . & . & .\\
1 & 1 & 1 & 1 & . & . & . & . & . & 1 & 1 & 1 & . & .\\
1 & . & . & . & 1 & . & . & 1 & . & 1 & . & . & 1 & .\\
1 & 1 & 1 & 1 & 1 & 1 & 1 & 1 & 1 & 1 & 1 & 1 & 1 & 1\\
\end{array}\right)
\end{equation*}
\caption{\label{e2p4}The matrix $M_{\EE_T,\Pp_T}$ for $n=4$.}
\end{figure}

\vfill
\newpage

~\vfill

\begin{figure}[ht]
\newcolumntype{C}{@{\hskip15pt}>{$}c<{$}}
\begin{equation*}
\left(\begin{array}{cC}
. & -1\\
1 & 1\\
\end{array}\right)
\qquad\qquad
\left(\begin{array}{cCCCC}
. & . & 1 & . & 1\\
. & . & -1 & . & .\\
. & . & . & -1 & -1\\
. & -1 & -1 & . & -1\\
1 & 1 & 1 & 1 & 1\\
\end{array}\right)
\end{equation*}
\caption{\label{e2h23}The matrices $M_{\EE_T,\HH_T}$ for $n=2,3$.}
\end{figure}

\vfill

\begin{figure}[ht]
\newcolumntype{C}{@{\hskip15pt}>{$}c<{$}}
\begin{equation*}
\left(\begin{array}{cCCCCCCCCCCCCC}
. & . & . & -1 & . & . & -1 & . & -1 & . & . & -1 & . & -1\\
. & . & . & 1 & . & . & . & . & . & . & . & 1 & . & .\\
. & . & . & . & . & . & 1 & . & . & . & . & . & . & .\\
. & . & . & . & . & . & . & 1 & 1 & . & . & . & 1 & 1\\
. & . & . & 1 & . & . & 1 & . & 1 & . & . & . & . & .\\
. & . & . & -1 & . & . & -1 & . & . & . & . & . & . & .\\
. & . & . & . & . & . & . & -1 & -1 & . & . & . & . & .\\
. & . & . & . & . & . & . & . & . & . & 1 & 1 & . & 1\\
. & . & . & . & . & . & . & . & . & -1 & -1 & -1 & -1 & -1\\
. & . & 1 & 1 & . & 1 & 1 & . & 1 & . & . & 1 & . & 1\\
. & . & -1 & -1 & . & . & . & . & . & . & . & -1 & . & .\\
. & . & . & . & -1 & -1 & -1 & -1 & -1 & . & . & . & -1 & -1\\
. & -1 & -1 & -1 & . & -1 & -1 & . & -1 & . & -1 & -1 & . & -1\\
1 & 1 & 1 & 1 & 1 & 1 & 1 & 1 & 1 & 1 & 1 & 1 & 1 & 1\\
\end{array}\right)
\end{equation*}
\caption{\label{e2h4}The matrix $M_{\EE_T,\HH_T}$ for $n=4$.}
\end{figure}

\vfill
\newpage

~\vfill

\begin{figure}[ht]
\newcolumntype{C}{@{\hskip15pt}>{$}c<{$}}
\begin{equation*}
\left(\begin{array}{cC}
1 & .\\
. & 1\\
\end{array}\right)
\qquad\qquad
\left(\begin{array}{cCCCC}
1 & . & . & . & .\\
. & 1 & . & . & .\\
. & . & . & 1 & .\\
. & . & 1 & 1 & .\\
. & . & . & . & 1\\
\end{array}\right)
\end{equation*}
\caption{\label{qp2q23}The matrices $M_{\Qq'_T,\Qq_T}$ for $n=2,3$.}
\end{figure}

\vfill

\begin{figure}[ht]
\newcolumntype{C}{@{\hskip15pt}>{$}c<{$}}
\begin{equation*}
\left(\begin{array}{cCCCCCCCCCCCCC}
1 & . & . & . & . & . & . & . & . & . & . & . & . & .\\
. & 1 & . & . & . & . & . & . & . & . & . & . & . & .\\
. & . & . & . & 1 & . & . & . & . & . & . & . & . & .\\
. & . & . & . & . & . & . & . & . & 1 & . & . & . & .\\
. & . & 1 & . & 1 & . & . & . & . & . & . & . & . & .\\
. & . & . & . & . & 1 & . & . & . & . & . & . & . & .\\
. & . & . & . & . & . & . & . & . & . & 1 & . & . & .\\
. & . & . & . & . & . & . & 1 & . & 1 & 1 & . & . & .\\
. & . & . & . & . & . & . & . & . & . & . & . & 1 & .\\
. & . & . & 1 & . & . & . & 1 & . & 1 & . & . & . & .\\
. & . & . & . & . & . & 1 & 1 & . & . & 1 & . & . & .\\
. & . & . & . & . & . & . & . & . & . & . & 1 & 1 & .\\
. & . & . & . & . & . & . & . & 1 & . & . & 1 & 1 & .\\
. & . & . & . & . & . & . & . & . & . & . & . & . & 1\\
\end{array}\right)
\end{equation*}
\caption{\label{qp2q4}The matrix $M_{\Qq'_T,\Qq_T}$ for $n=4$.}
\end{figure}

\vfill
\newpage

~\vfill

\begin{figure}[ht]
\newcolumntype{C}{@{\hskip15pt}>{$}c<{$}}
\begin{equation*}
\left(\begin{array}{cC}
1 & .\\
-1 & 1\\
\end{array}\right)
\qquad\qquad
\left(\begin{array}{cCCCC}
1 & . & . & . & .\\
-1 & 1 & . & . & .\\
. & -1 & 1 & . & .\\
-1 & . & . & 1 & .\\
1 & . & -1 & -1 & 1\\
\end{array}\right)
\end{equation*}
\caption{\label{m2q23}The matrices $M_{\Prima_T,\Qq_T}$ for $n=2,3$.}
\end{figure}

\vfill
\begin{figure}[ht]
\newcolumntype{C}{@{\hskip15pt}>{$}c<{$}}
\begin{equation*}
\left(\begin{array}{cCCCCCCCCCCCCC}
1 & . & . & . & . & . & . & . & . & . & . & . & . & .\\
-1 & 1 & . & . & . & . & . & . & . & . & . & . & . & .\\
. & -1 & 1 & . & . & . & . & . & . & . & . & . & . & .\\
. & . & -1 & 1 & . & . & . & . & . & . & . & . & . & .\\
-1 & . & . & . & 1 & . & . & . & . & . & . & . & . & .\\
1 & . & -1 & . & -1 & 1 & . & . & . & . & . & . & . & .\\
. & . & 1 & -1 & . & -1 & 1 & . & . & . & . & . & . & .\\
. & . & . & . & -1 & . & . & 1 & . & . & . & . & . & .\\
. & . & . & . & 1 & . & -1 & -1 & 1 & . & . & . & . & .\\
-1 & . & . & . & . & . & . & . & . & 1 & . & . & . & .\\
1 & -1 & . & . & . & . & . & . & . & -1 & 1 & . & . & .\\
. & 1 & . & -1 & . & . & . & . & . & . & -1 & 1 & . & .\\
1 & . & . & . & . & . & . & -1 & . & -1 & . & . & 1 & .\\
-1 & . & . & 1 & . & . & . & 1 & -1 & 1 & . & -1 & -1 & 1\\
\end{array}\right)
\end{equation*}
\caption{\label{m2q4}The matrix $M_{\Prima_T,\Qq_T}$ for $n=4$.}
\end{figure}

\vfill
\newpage

~\vfill

\begin{figure}[ht]
\newcolumntype{C}{@{\hskip15pt}>{$}c<{$}}
\begin{equation*}
\left(\begin{array}{cC}
1 & -1\\
. & 1\\
\end{array}\right)
\qquad\qquad
\left(\begin{array}{cCCCC}
1 & -1 & . & -1 & 1\\
. & 1 & -1 & . & .\\
. & . & 1 & . & -1\\
. & . & . & 1 & -1\\
. & . & . & . & 1\\
\end{array}\right)
\end{equation*}
\caption{\label{fm2q23}The matrices $M_{\Primb_T,\Qq_T}$ for $n=2,3$.}
\end{figure}

\vfill

\begin{figure}[ht]
\newcolumntype{C}{@{\hskip15pt}>{$}c<{$}}
\begin{equation*}
\left(\begin{array}{cCCCCCCCCCCCCC}
1 & -1 & . & . & -1 & 1 & . & . & . & -1 & 1 & . & 1 & -1\\
. & 1 & -1 & . & . & . & . & . & . & . & -1 & 1 & . & .\\
. & . & 1 & -1 & . & -1 & 1 & . & . & . & . & . & . & .\\
. & . & . & 1 & . & . & -1 & . & . & . & . & -1 & . & 1\\
. & . & . & . & 1 & -1 & . & -1 & 1 & . & . & . & . & .\\
. & . & . & . & . & 1 & -1 & . & . & . & . & . & . & .\\
. & . & . & . & . & . & 1 & . & -1 & . & . & . & . & .\\
. & . & . & . & . & . & . & 1 & -1 & . & . & . & -1 & 1\\
. & . & . & . & . & . & . & . & 1 & . & . & . & . & -1\\
. & . & . & . & . & . & . & . & . & 1 & -1 & . & -1 & 1\\
. & . & . & . & . & . & . & . & . & . & 1 & -1 & . & .\\
. & . & . & . & . & . & . & . & . & . & . & 1 & . & -1\\
. & . & . & . & . & . & . & . & . & . & . & . & 1 & -1\\
. & . & . & . & . & . & . & . & . & . & . & . & . & 1\\
\end{array}\right)
\end{equation*}
\caption{\label{fm2q4}The matrix $M_{\Primb_T,\Qq_T}$ for $n=4$.}
\end{figure}

\vfill
\newpage

%
%

\end{document}